\documentclass[12pt,a4paper]{amsart}
\usepackage{amsmath,amsthm,amsfonts,amssymb,mathrsfs}
\usepackage{subfigure} 
\usepackage{graphicx}
\usepackage{color}
\newcommand\nofig[1]{#1}

\newcommand\shorter[1]{}

\newcommand\R{\mathbb{R}}
\newcommand\C{\mathbb{C}}
\newcommand\N{\mathbb{N}}
\newcommand\A{\mathcal{A}}
\newcommand\I{\mathcal{I}}
\newcommand\F{\mathcal{F}}

\renewcommand\H{\mathcal{H}}
\newcommand\wH{\overline{H}}
\newcommand\wT{\widetilde{T}}

\newcommand\U{\mathcal{U}}
\newcommand\V{\mathcal{V}}

\renewcommand\L{\mathcal{L}}
\renewcommand\Re{{\rm Re \,}}
\def\u{\bar{u}}
\def\v{\bar{v}}
\def\w{\bar{w}}

\newtheorem{theorem}{Theorem}
\newtheorem{cor}[theorem]{Corollary}

\newtheorem{defi}[theorem]{Definition}

\newtheorem{lemma}[theorem]{Lemma}
\newtheorem{prop}[theorem]{Proposition}

\numberwithin{equation}{section}
\numberwithin{theorem}{section}
\numberwithin{figure}{section}

\theoremstyle{remark}
\newtheorem{rem}[theorem]{Remark}

\usepackage{delarray}
\usepackage{enumerate}

\begin{document}

\title[Reaction-diffusion model of early carcinogenesis]{Unstable patterns 
in reaction-diffusion model  \\ of early carcinogenesis }

\author[A.~Marciniak-Czochra]{Anna Marciniak-Czochra}
\address[A. Marciniak-Czochra]{
 Interdisciplinary Center for Scientific Computing (IWR), Institute of Applied Mathematics and BIOQUANT, University of Heidelberg, 69120 Heidelberg, Germany}
\email{anna.marciniak@iwr.uni-heidelberg.de}
\urladdr{http://www.iwr.uni-heidelberg.de/groups/amj/People/Anna.Marciniak}

\author[G.~Karch]{Grzegorz Karch}
\address[G. Karch]{
 Instytut Matematyczny, Uniwersytet Wroc\l awski,
 pl. Grunwaldzki 2/4, 50-384 Wroc\-\l aw, Poland}
\email{grzegorz.karch@math.uni.wroc.pl}
\urladdr{http://www.math.uni.wroc.pl/~karch}

\author[K.~Suzuki]{Kanako Suzuki}
\address[K. Suzuki]{
Graduate School of Information Sciences, Tohoku University,
6-3-09 Aramaki-aza-Aoba, Aoba-ku, Sendai 980-8579, Japan}
\email{kasuzu@m.tohoku.ac.jp}

\date{\today}


\begin{abstract}
Motivated by numerical simulations showing the emergence of either periodic or irregular patterns, we explore a mechanism of pattern formation arising in the processes described by a system of a single reaction-diffusion equation coupled with ordinary differential equations. We focus on a basic model of early cancerogenesis proposed by Marciniak-Czochra and Kimmel [Comput.~Math.~Methods Med.~{\bf 7} (2006), 189--213],
[Math.~Models Methods Appl.~Sci.~{\bf 17} (2007), suppl., 1693--1719],
but the theory we develop applies to a wider class of pattern formation models with an autocatalytic non-diffusing component. The model exhibits diffusion-driven instability (Turing-type instability). However, we prove that all Turing-type patterns, {\it i.e.,} regular stationary solutions,  are unstable
in the Lyapunov sense. Furthermore, we show existence of discontinuous stationary solutions, which are also unstable. 
\end{abstract}
\maketitle

\section{Introduction}

In this paper, we explore a mechanism of pattern formation arising in  processes described by a system of a single reaction-diffusion equation coupled with ordinary differential equations. Such systems of equations arise from modeling of interactions between cellular processes and diffusing growth factors.

 A rigorous derivation, using methods of asymptotic analysis (homogenization), of the macroscopic reaction-diffusion models
describing the interplay between the nonhomogeneous cellular
dynamics and the signaling molecules diffusing in the intercellular
space has been recently published in \cite{MCP}. It was shown that 
receptor-ligand binding processes can be modeled by reaction-diffusion equations coupled with ordinary differential equations in the case when all membrane processes are homogeneous within 
the membrane, which seems to be the case in most of processes. 
More precisely,  such  {\it receptor-based} models
can be represented by the following initial-boundary value 
problem
$$
\begin{aligned}
u_t & =   f(u,v),&&\\
v_t & =   D \Delta v+g(u,v) &\quad &\text{in}\quad \Omega,\\
\partial_{n}v & =  0  &&\text{on}\quad \partial\Omega,\\
u(x,0) & =   u_{0}(x),&&\\
v(x,0) & =  v_{0}(x),&&
\end{aligned}
$$ 
where $u$ and $v$  are vectors  of  variables,  $D$ is a diagonal matrix with positive coefficients on the diagonal,  the symbol
$\partial_{n}$ denotes the normal derivative (no-flux condition),
and
$\Omega$ is a bounded region.
It was shown that, if  homogeneity of the processes on the membrane does not hold, equations with additional integral terms  are obtained, see \cite{MCP}.

One of  possible mechanisms of pattern formation in such models is {\it the diffusion-driven instability} (DDI), called also {\it the Turing-type instability}. Let us recall that the diffusion-driven instability arises in a reaction-diffusion system, when there exists a spatially homogeneous solution, which is asymptotically stable  with respect to spatially homogeneous perturbations, but  unstable to spatially inhomogeneous perturbation.  Models with DDI describe the process of destabilization of a stationary spatially homogeneous state and  evolution of spatially heterogeneous structures, which converge to a spatially heterogeneous steady state. The systems  with DDI were prevalent  in the modeling literature since the seminal paper of  Turing \cite{Turing} and  have provided explanations of pattern formation in a variety of biological systems, see {\it e.g.} \cite{Murray} and references therein.

The existing qualitative theory of such systems is mainly focused on stability conditions for homogeneous steady states. Moreover, the majority of theoretical studies focus on the analysis of the non-degenerated reaction-diffusion systems, {\it i.e.}~with a strictly positive  diffusion coefficient in each equation.  However, in many biological applications, it is relevant to consider receptor-based systems involving ordinary differential equations.  An interesting class of such systems consists of only a single reaction-diffusion equation coupled with a system of ordinary differential equations.  The existence and stability of spatially heterogeneous patterns arising in models exhibiting diffusion-driven
instability, but consisting of only one reaction-diffusion equation, is a mathematically interesting issue.
Such models are very different from classical Turing-type models and the spatial structure of the pattern emerging from the destabilization of the spatially homogeneous steady state cannot be concluded based on linear stability analysis \cite{Murray}. 
Asymptotic analysis of systems with such a dispersion relation seems to be an open problem. 

In particular, such systems arise in the modeling of the growth of a spatially-distributed cell population, which proliferation is controlled by endogenous or exogenous growth factors diffusing in the extracellular medium and binding to cell surface as proposed by the first author  and Kimmel in the series of recent papers \cite{MCK06, MCK07, MCK08}. The generic model has  the following mathematical form of mixed reaction-diffusion and ordinary differential equations
\begin{align}
u_t&=\left(a_p(u,v) -d_c\right) u, \label{eq1-intro}\\
v_t &=\alpha(u) w -d v  -d_b v, \label{eq2-intro}\\
w_t &= D \Delta w  -d_g w-\alpha(u) w  +d v +\kappa(u). \label{eq3-intro}
\end{align}
Here, unknown functions $u(x,t)$, $v(x,t)$, $w(x,t)$ describe the densities of cells, free and bound growth factor molecules, respectively,  distributed over a certain bounded domain.  
In system  \eqref{eq1-intro}-\eqref{eq3-intro}, 
the proliferation rate $a_p(u,v)$,  the growth factor binding rate $\alpha(u)$, and the production of the free 
growth factor by cells $\kappa(u)$
are given functions. We use the subscripts in the parameters $ d_c$, $d_b$, $d_g$ in order to emphasise that they denote the degradation  factors of cells, bound molecules and growth molecules, respectively.

System \eqref{eq1-intro}-\eqref{eq3-intro} is usually  supplemented 
with the homogeneous Neumann (zero  flux) boundary conditions for the function $w=w(x,t)$
and with nonnegative initial conditions.

 Based on preliminary  mathematical analysis and using numerical simulations, 
the authors of \cite{MCK06, MCK07, MCK08} found conditions 
for the existence of a positive spatially homogeneous steady state exhibiting diffusion-driven instability. In numerical simulations instability of the constant steady state  leads to the emergence of growth patterns concentrated
around discrete points along the spatial coordinate, which take the mathematical form of
spike-type spatially inhomogeneous solutions. This multifocality is as expected from the
field theory of carcinogenesis. Numerical simulations showed qualitatively new patterns of behavior of solutions, including, in some cases, a strong dependence of the emerging pattern on initial conditions and quasi-stability followed by rapid growth of solutions.

The main goal of this work is to develop mathematical theory that allows us to study properties of
solutions to this kind of initial-boundary value problems.
In this work,  for simplicity of the exposition, we assume that cells occupy the interval $x\in [0,L]$ and  we consider 
a particular form of functions   $a_p$, $\alpha$, $\kappa$, see the next section. 
Our results, however, hold true in a much more general case what we have systematically emphasized 
in our construction of patterns (see Section \ref{sec:pattern}) and in the proof of their instability (Section \ref{sect:spect}).

The paper is organized as follows. In Section \ref{Sec2}, main results of this paper are formulated. Section \ref{sec:existence} is devoted to the existence of  solutions. Section \ref{sec:large} provides 
preliminary results on the large time behavior of solutions and, in particular,
a criterion on the pointwise extinction of solutions. In Section \ref{sec:pattern}, nonhomogeneous 
 stationary solutions are characterized and, in  Section \ref{sect:spect},  instability of all stationary solutions  is shown using linear stability analysis. The paper is supplemented by  Appendix devoted to the analysis of the 
corresponding
kinetics system (system of ordinary differential equations).

\subsection*{Notation}
The usual norm of the Lebesgue space 
$L^p (0,1)$ is denoted by $\|\cdot\|_p$ for any $p \in [1,\infty]$ and $W^{k,p}(0,1)$ is the corresponding Sobolev space. 
The constants (always independent
of $x$ and $t $) will be
denoted by the
same letter $C$, even if they may vary from line to line.
Sometimes, we write,  {\it e.g.},  $C=C(p,q,r, ...)$ when we want to
emphasize
the dependence of $C$ on parameters~$p,q,r, ...$.


\section{Results and coments}\label{Sec2}

\subsection{Statement of the problem}
In this work, we consider system \eqref{eq1-intro}-\eqref{eq3-intro} on 
 the bounded interval $x\in [0,L]$. Moreover,  in our mathematical analysis,
we assume that the proliferation rate  has the Hill function form 
$$
a_p(u,v) =a \frac{v/{u}}{1+{v}/{u}}= \frac{av}{u+v}
$$
 with a given constant $a>0$. We consider the quadratic  growth factor binding rate 
$\alpha(u)=u^2$, which follows from conditions for the diffusion driven instability,
see \cite{MCK06} and the next subsection
for more details. Moreover, we assume that the production of the free growth factor is constant,
$\kappa(u)\equiv \kappa_0\geq 0$. This assumption is imposed for simplicity of the exposition; note that a more general case is also interesting
from the modeling point of view,  see  \cite{MCK07}.
Finally, we introduce the diffusion
coefficient $1/\gamma$, which is a composite
parameter including the diffusion constant $D>0$ and the scaling parameter $L>0$, namely, we set
$\gamma =L^2/D$.

To summarize, we  study the following system of three ordinary/partial differential equations 
\begin{align}
u_t&=\Big(\frac{a v}{u+v} -d_c\Big) u \qquad \text{for}\ x\in [0,1], \; t>0, \label{eq1}\\
v_t &=-d_b v +u^2 w -d v \qquad  \text{for}\ x\in [0,1], \; t>0, \label{eq2}\\
w_t &= \frac{1}{\gamma} w_{xx} -d_g w -u^2 w +d v +\kappa_0 \qquad \text{for}\ x\in (0,1), \; t>0 \label{eq3}
\end{align}
supplemented with the homogeneous Neumann (zero  flux) boundary conditions for the function $w=w(x,t)$
\begin{equation}\label{N}
w_x(0,t)=w_x(1,t)=0 \quad \mbox{for all} \quad t>0
\end{equation}
and with nonnegative initial conditions
\begin{equation}
u(x,0)=u_0(x), \quad v(x,0)=v_0(x), \quad w(x,0)=w_0(x).\label{ini}
\end{equation}
Here, the letters $a , d_c,d_b, d_g, d, \gamma, \kappa_0 $  denote positive constants.

In Section \ref{sec:existence}, we  show that the initial-boundary  value problem \eqref{eq1}--\eqref{ini}
  has  a unique and global-in-time solution
for a large class of nonnegative initial conditions.  
Such results on the global-in-time existence, the regularity of solutions, and their positivity for all $t>0$  
are rather standard for reaction-diffusion equations 
with non-zero diffusion in each equation, see {\it e.g.} 
\cite{pierre}
and the references therein. A more careful analysis is required 
in the case of the ODE-PDE system, hence, 
we state these results for the completeness of the exposition. 

\subsection{Diffusion-driven instability}
We are interested in systems with the diffusion-driven instability. As it was stated in  \cite[Prop.~3.1]{MCK08}, a generic system of two ordinary differential equations  coupled with a reaction-diffusion equation, and such that   $a_{ii} < 0$ for $i = 1, 2, 3$ and  $a_{12}a_{21} > 0$,  exhibits the diffusion-driven instability if there exists a positive, spatially constant steady state, for which the following conditions are satisfied
\begin{eqnarray}
-{\rm tr}\,(\A)>0,\label{cond1}&&\\
-{\rm tr}\,(\A)\sum_{i<j}\det (\A_{ij})+\det(\A)>0,\label{cond2}&&\\
-\det(\A)>0,\label{cond3}&&\\
-\det(\A_{12})>0,\label{cond4}&&
\end{eqnarray}
where $\A=\big(a_{ij}\big)_{i,j=1,2,3}$  is the Jacobian matrix of the system without diffusion linearized around this
constant  positive equilibrium and  $\A_{ij}$ is a submatrix of $\A$ consisting of the $i$-th and $j$-th column and $i$-th and $j$-th row.

Conditions \eqref{cond1}-\eqref{cond3} are necessary for the stability of the steady state in the absence
of diffusion.  Inequality \eqref{cond4} is a sufficient and necessary condition for destabilization of
this steady state in every system with  all $a_{ii} < 0$ for $i = 1, 2, 3$ and  $a_{12}a_{21} > 0$ (for the proof see \cite{MCK06}).
These conditions guarantee that the model \eqref{eq1-intro}-\eqref{eq3-intro} exhibits diffusion-driven instability if the function $\alpha(u)$ evaluated at the steady state $\u$ satisfies  $\alpha(\bar u) - \bar u \alpha'(\bar u) <0$, what, in particularly, always holds for $\alpha(u) =u^2$ in model \eqref{eq1}-\eqref{eq3}.

In the particular case of system \eqref{eq1}-\eqref{eq3}, conditions \eqref{cond1}-\eqref{cond4} hold true 
for the constant steady state $(\u_{-},\v_{-},\w_{-})$ defined in \eqref{const-pm} (see the next subsection) and for the matrix $\A$ stated in \eqref{A0} with $W=\w_{-}$.
Note that condition \eqref{cond4} leads to 2 eigenvalues of the matrix $\A_{12}$ of an opposite sign.
The positive eigenvalue of $\A_{12}$ plays the crucial role in the proof of instability of all
steady states of system  \eqref{eq1}-\eqref{eq3}, see Subsection \ref{subsec:2.5}.

\subsection{Preliminary properties of solutions}
We begin our study of qualitative properties of solutions to problem \eqref{eq1}--\eqref{ini} 
by considering  $x$-independent  solutions which
satisfy  the corresponding {\it kinetic system} of the three ordinary
differential equations (see \eqref{heq1}--\eqref{heq3}, below). %
Let us briefly summarize our results on the kinetic system that  we prove  in Appendix~A.

\begin{itemize}
\item For every constant initial condition $\u(0)\geq 0$, $\v(0)\geq 0$, and $\w(0)\geq 0$, the problem \eqref{eq1}--\eqref{ini}
has a unique $x$-independent global-in-time solution $(\u(t), \v(t), \w(t))$. This solution 
stays, for $t\geq 0$, in a bounded set, see Proposition \ref{app-pro1}.
Here,  if $\u(t_0)=\v(t_0)=0$ for  certain $t_0\geq 0$, 
the right-hand side of equation \eqref{eq1} is satisfied  in the limit sense, namely, when $\v\searrow 0$ and $\u\searrow 0$.

\item
The constant vector 
\begin{equation}\label{trivial:steady}
(\u_0,\v_0,\w_0)\equiv  \Big(0, 0, \frac{\kappa_0}{d_g}\Big)
\end{equation}
is a {\it trivial  steady state} of system  \eqref{eq1}--\eqref{eq3}. %
This is an asymptotically stable solution not only of the kinetic system \eqref{heq1}--\eqref{heq3}, but also of the reaction-diffusion equations 
\eqref{eq1}-\eqref{N},
see Corollary \ref{cor:stab:trivial} and Proposition  \ref{prop:a<dc}.
Moreover,  in Theorems \ref{app-th1}--\ref{app-th4} of Appendix, we describe  convergence rates
of solutions to the kinetic system towards the trivial steady states.

\item
Assume that $a>d_c$ and 
\begin{equation}\label{triangle}
 \kappa_0^2\geq \Theta, \qquad \text{where}\quad \Theta= 4d_g 
d_b
\frac{d_c^2 (d_b+d)}{(a -d_c)^2}.
\end{equation}
Then, system \eqref{eq1}--\eqref{eq3} has constant positive stationary solutions 
$(\u_{\pm},\v_{\pm},\w_{\pm}),$
where 
\begin{equation}\label{const-pm}
\w_{\pm}= \frac{\kappa_0\pm \sqrt{\kappa_0^2-\Theta }}{2 d_g}, \quad 
\v_\pm=\frac{d_c^2 (d_b+d)}{(a -d_c)^2}\; \frac{1}{\w_\pm}, \quad 
\u_\pm=\frac{a -d_c}{d_c}\; \v_\pm.
\end{equation}

\item The stationary solution $(\u_{-},\v_{-},\w_{-})$ of the kinetic system 
\eqref{heq1}--\eqref{heq3} is  asymptotically  stable, namely, the linearization matrix of this system at the steady state 
$(\u_{-},\v_{-},\w_{-})$ has all eigenvalues with negative real parts. 
On the other hand, the steady state 
$(\u_{+},\v_{+},\w_{+})$ is an unstable solution of \eqref{heq1}--\eqref{heq3}, because
its linearization matrix has one positive eigenvalue, see Theorem \ref{app-lem4} and Corollary \ref{app-cor1}, for more details.

\end{itemize}

From now on, we consider solutions of problem \eqref{eq1}--\eqref{ini}
that are not necessarily  space-inhomogeneous. %
Unlike in the case of the kinetic system (see Proposition \ref{app-pro1} in Appendix), 
in view of numerical simulations, 
we do not expect  that such a  solution  
stays in a certain invariant region as
it is discussed in the monograph by Smoller \cite[Ch. 14, \S B]{Smoller}.
However, in the following theorem, we prove  that integrals of all nonnegative solutions of 
\eqref{eq1}--\eqref{ini} enter, as $t\to\infty$, into a certain invariant set, which is independent of 
initial conditions.
Notice that  the mass of a solution 
is controlled for all $t>0$  as in the survey article \cite{pierre},
were systems with non-zero diffusion in each equation were discussed.

\begin{theorem}\label{thm:inv}
Assume that $(u, v, w)$ is a nonnegative 
global-in-time solution of problem  \eqref{eq1}-\eqref{ini}  corresponding to
a bounded initial condition $(u_0,v_0,w_0)\in L^1(0,1)\times L^1(0,1)\times L^\infty(0,1)$. 
Denote $\mu=\min\{d_g,d_b\}>0$. Then, the following estimates hold
\begin{align}
\limsup_{t\to\infty} \int_0^1 u(x,t)\;dx&\leq  \min\left\{\frac{\kappa_0}{\mu },  \frac{a \kappa_0}{d_c\mu }\right\}, 
 \label{limsup:1}\\
\limsup_{t\to\infty} \int_0^1 v(x,t)\;dx&\leq \frac{\kappa_0}{\mu},\label{limsup:2}\\
\limsup_{t\to\infty} \|w(t)\|_\infty &\leq {\kappa_0}\left(\frac{Cd}{\mu d_g^{1/2}}+\frac{1}{d_g}\right),\label{limsup:3}
\end{align}
with a numeric constant $C>0$ independent of problem \eqref{eq1}-\eqref{ini}.
\end{theorem}

This theorem is proved at the end of Section  \ref{sec:existence}.
Notice that it holds also true for $\kappa_0 = 0$, hence, in this particular case, each %
nonnegative solution of problem \eqref{eq1}-\eqref{ini}
satisfies 
$\big(\|u(t)\|_1,\|v(t)\|_1,\|w(t)\|_\infty \big)\to (0,0,0)
$
as  $t\to \infty$.

\medskip

Next, we discuss the stability of the trivial steady state \eqref{trivial:steady} as a solution of the reaction-diffusion equations
 \eqref{eq1}-\eqref{N}.
First of all, we have to emphasize that, 
under the assumption 
$a<d_c$,
 each nonnegative solution of problem \eqref{eq1}-\eqref{ini} converges exponentially
towards  $(0,0, \kappa_0/d_g)$, see Proposition \ref{prop:a<dc} below.
Hence, in this work,  we have  always to assume that $a\geq d_c$  to observe a nontrivial large time behavior of solutions.

It is clear from equation \eqref{eq1} that if $u_0(x)=0$ for some $x\in [0,1]$ than $u(x,t)=0$ for all
$t\geq 0$ and, by equation \eqref{eq2}, we have $v(x,t)\to 0$ as $t\to \infty$.
Furthermore, below in Proposition~\ref{prop:zero}, we show that, for each $x\in [0,1]$, 
the following two conditions 
\begin{equation}\label{uv0}
\lim_{ t\to\infty} u(x,t) =0 
\qquad \text{and} \qquad  \lim_{ t\to\infty} v(x,t)=0
\end{equation}
are equivalent. The following  result plays a fundamental role in  understanding of the pattern formation
described by problem \eqref{eq1}-\eqref{ini} 
(see Remark \ref{rem:disc}, below)
and says that relations \eqref{uv0} hold true if $u_0(x)$ and $v_0(x)$ are not
too large.

\begin{theorem}\label{thm:stab:trivial}
Assume that $a>d_c$ and other coefficients in system \eqref{eq1}-\eqref{ini} are positive and arbitrary.
Fix $x \in [0,1]$ and assume that there exist constants  $K_w>0$ and $M>0$
satisfying
\begin{equation} \label{as:M}
MK_w \left(1+\frac{d_c}{a}\right)^2 \leq (d_b+d) \left(\frac{d_c}{a}\right)^2.
\end{equation}
and  such that for all $t>0$
\begin{equation}\label{u:M}
0\leq w(x, t) \leq  K_w  \quad  \text{and} \quad 0\leq u_0(x) \leq M, \quad  0\leq v_0(x)<\left(\frac{d_c}{a}\right)^2 M.
\end{equation} 
Then, $(u(x,t),v(x,t))\to (0,0)$ as $t\to \infty$. This is the uniform convergence for all $x\in [0,1]$, for which 
inequalities \eqref{u:M} are satisfied.
\end{theorem}

\begin{rem}
It follows from Theorem \ref{thm:inv} that $\|w(t)\|_\infty$ is  a bounded function of $t>0$, hence, one can always find a constant $K_w>0$
with the property stated in Theorem \ref{thm:stab:trivial}. %
In fact, this constant can be chosen explicitly, see Remark \ref{bound-w-rem} at the
 end of Section \ref{sec:existence}. %
\end{rem}

We conclude this subsection by a corollary of Theorem  \ref{thm:stab:trivial} on the asymptotic stability of the trivial steady state.

\begin{cor}\label{cor:stab:trivial}
Let $a>d_c$ and other coefficients in system \eqref{eq1}-\eqref{ini} be positive.
Assume that there exist constants  $K_w>0$ and $M>0$ satisfying \eqref{as:M} and \eqref{u:M}  for
all $x\in [0,1]$.
Then $\big(u(x,t), v(x,t), w(x,t)\big)\to (0,0, \kappa_0/d_g)$ as $t\to\infty$ uniformly in $x$.
\end{cor}

Here, it has to be emphasized that system  \eqref{eq1}-\eqref{ini} cannot be linearized around the
solution $(0,0, \kappa_0/d_c)$, because of the singularity at $u=v=0$ in equation \eqref{eq1}.
Hence, to prove the convergence of a solution towards the trivial  steady state
a direct  method
(based on the ordinary differential equations \eqref{eq1}-\eqref{eq2})
 was invented, see  the proof of Theorem  \ref{thm:stab:trivial} in Section \ref{sec:large}.

\subsection{Existence and nonexistence  of patterns}

Here, 
 we describe  stationary solutions of  \eqref{eq1}-\eqref{N}, namely, 
 functions  $(U(x), V(x), W(x))$ that satisfy the system
\begin{align}
\left(\frac{a  V}{U+V} -d_c\right) U&=0, \label{seq1}\\
-d_b V +U^2 W -d V&=0,\label{seq2}\\
 \frac{1}{\gamma} W_{xx} -d_g W -U^2 W +d V +\kappa_0&=0 \label{seq3}
\end{align}
and the boundary condition
\begin{equation}\label{sN}
W_x(0)=W_x(1)=0. 
\end{equation}

Let us recall that, for every $\gamma>0$, system \eqref{seq1}-\eqref{sN} has the trivial solution %
$(\u_0,\v_0,\w_0)=(0, 0, {\kappa_0}/{d_g})$. Moreover, if $\kappa_0^2\geq \Theta$, see \eqref{triangle},
we have two other constant stationary solutions $(\u_{\pm},\v_{\pm},\w_{\pm})$  defined in \eqref{const-pm}.

In the nonhomogeneous case, first, we consider 
$U(x)$ and $V(x)$, which  are positive for all $x\in [0,1]$.
Hence, 
under the assumption $a>d_c$,  using equations \eqref{seq1}-\eqref{seq2} we obtain the following identities
\begin{equation}\label{s-uv}
U(x)=\frac{a -d_c}{d_c}\; V(x)  \quad \mbox{and}\quad V(x)=\frac{d_c^2 (d_b+d)}{(a -d_c)^2}\; \frac{1}{W(x)}.
\end{equation}
Moreover, adding equations \eqref{seq2} and \eqref{seq3}, and using relations \eqref{s-uv} we obtain
the following boundary-value problem for the function $W$:
\begin{align}
\frac{1}{\gamma} W'' - d_gW - d_b
\frac{d_c^2 (d_b+d)}{(a -d_c)^2}\; \frac{1}{W}
 +\kappa_0&=0,\label{sw1}\\
 W'(0)=W'(1)&=0.\label{sw2} 
\end{align}
Here, to simplify our notation, we introduce the function
\begin{equation}\label{h}
h(w)=- d_gw - d_b
\frac{d_c^2 (d_b+d)}{(a -d_c)^2}\; \frac{1}{w}
 +\kappa_0
\end{equation}
and we define number
\begin{equation}\label{gamma0}
\gamma_0\equiv \frac{\pi^2}{h'(\w_{-})}>0,
\end{equation}
where $\w_{-}$ is the constant from \eqref{const-pm}.

First, let us  summarize our results (proved in Section \ref{sec:pattern}) on the {\it nonexistence} of solutions either to system  
\eqref{seq1}-\eqref{sN} or to boundary value problem \eqref{sw1}-\eqref{sw2}.

\begin{itemize}
\item Under the assumption $a\leq d_c$, the trivial solution $(0, 0, {\kappa_0}/{d_g})$ is the only nonnegative solution 
of problem \eqref{seq1}-\eqref{sN}, see Proposition~\ref{prop:stat}.

\item If $a> d_c$ and $\kappa_0^2<\Theta$, then  the boundary value problem \eqref{sw1}-\eqref{sw2} has no positive solutions,
see Proposition~\ref{thm:stat-1}.i. Hence again, the trivial solution is the only nonnegative solution 
of  \eqref{seq1}-\eqref{sN}.

\item
If $a> d_c$ and $\kappa_0^2=\Theta$,
then the constant function $w\equiv \kappa_0/d_c$ is the only solution
of   \eqref{sw1}-\eqref{sw2}, see Proposition~\ref{thm:stat-1}.ii.

\item Let $a> d_c$ and $\kappa_0^2>\Theta$.
The parameter $\gamma_0$ is the critical value of the diffusion coefficient in the sense that
for $\gamma\leq \gamma_0$,  the boundary value problem \eqref{sw1}-\eqref{sw2} has no non-constant solutions., see
Theorem~\ref{thm:nopattern}. 

\end{itemize}

Hence, continuous patterns, {\it i.e.,~non-constant continuous positive solutions}, of 
 system  
\eqref{seq1}-\eqref{sN} can exist only if $a> d_c$, $\kappa_0^2>\Theta$, and $\gamma>\gamma_0$. 
Here, we show that this is indeed the case.

\begin{defi}
Let $k\in \N$ and $k\geq 2$. We call a function $W\in C([0,1])$ a periodic function on $[0,1]$ with $k$ modes 
if
$W=W(x)$ is monotone on $\left[0, \frac{1}{k}\right]$ and if
\begin{equation}\label{modes}
W(x) =
\left\{
\begin{array}{ccc}
W\left(x-\frac{2j}{k}\right)& \; \text{for}\;  & x \in \left[\frac{2j}{k}, \frac{2j+1}{k}\right]\\[3pt]
W\left(\frac{2j+2}{k}-x \right)& \; \text{for}\;  & x \in \left[\frac{2j+1}{k}, \frac{2j+2}{k}\right]
\end{array}
\right.
\end{equation}
for every $j\in \{0,1,2,3,...\}$ such that $2j+2\leq k$.
\end{defi}

We are now in a position to describe all possible nonnegative solution of 
the boundary value problem \eqref{sw1}-\eqref{sw2}.

\begin{theorem}[Existence of continuous patterns] \label{thm:stationary}
Assume that  $a>d_c$,  $\kappa_0^2>\Theta$,  and
  $\gamma>\gamma_0$, where $\gamma_0$ is defined in \eqref{gamma0}.
Consider the biggest $n\in\N$ such that 
$\gamma>n^2\gamma_0$.  
Then,  the boundary value problem \eqref{sw1}-\eqref{sw2}  has the following  solutions:
\begin{itemize}

\item the constant steady states $ \w_\pm$,

\item a unique strictly increasing solution and a unique strictly decreasing solution,

\item for each $k \in \{2,...,n\}$,
a unique periodic solution $W_k$ with $k$ modes
that is increasing on $[0, \frac{1}{k}]$
as well as 
 its symmetric counterpart: $\widetilde W_k(x)\equiv W_k(1-x) $.

\end{itemize}
There are no other positive solutions of  the boundary value problem \eqref{sw1}-\eqref{sw2}.
\end{theorem}

All positive solutions of problem \eqref{sw1}-\eqref{sw2} with 
with $\gamma$ satisfying assumptions of Theorem \ref{thm:stationary} for 
$n=3$   are sketched in Fig. \ref{fig:0}.

%
\nofig{

\begin{figure}
  \setlength{\unitlength}{1mm}
{\footnotesize
\begin{picture}(50,75)
\put(0,50){\includegraphics[width=20mm]{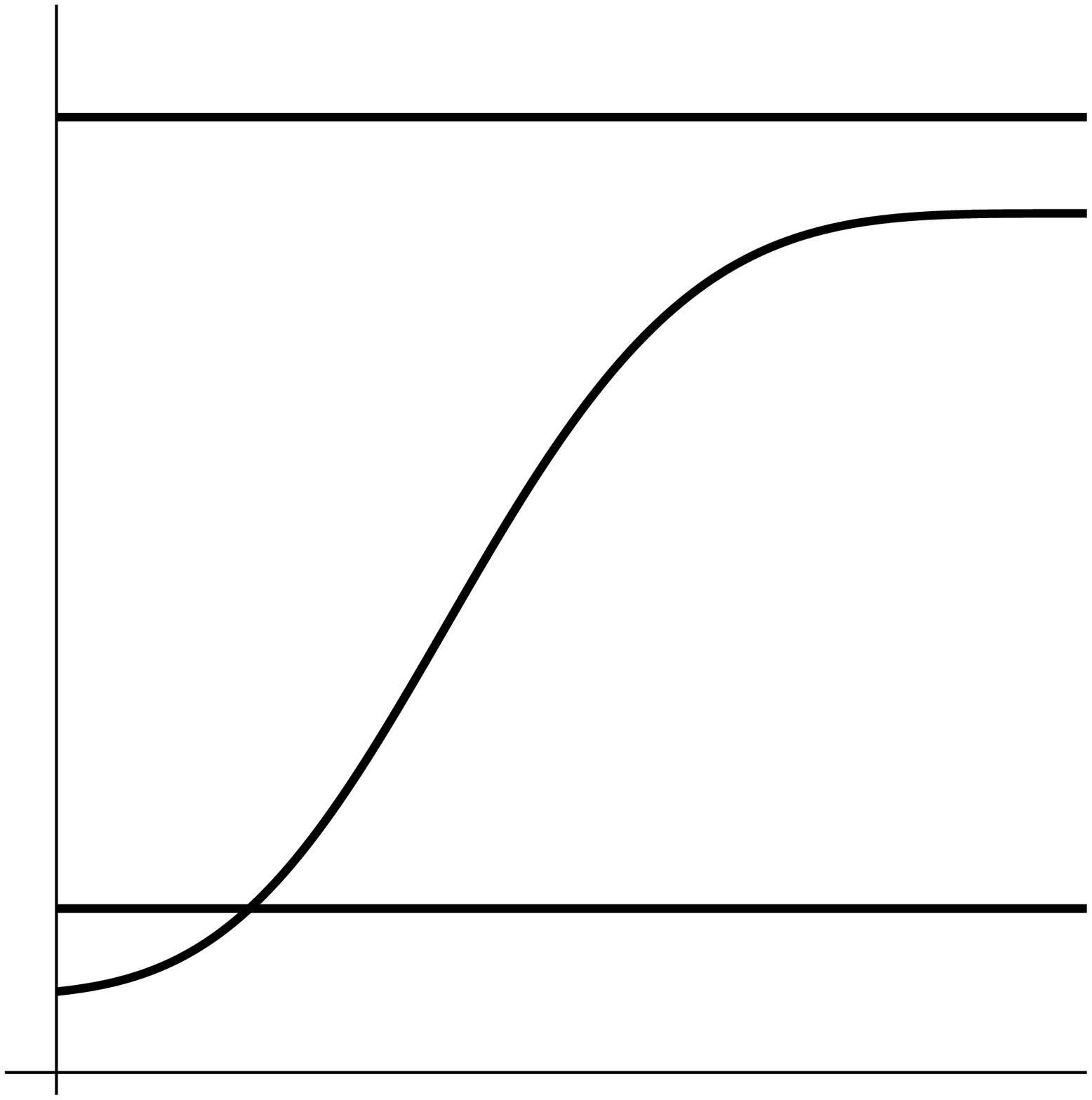}}
\put(26,50){\includegraphics[width=20mm]{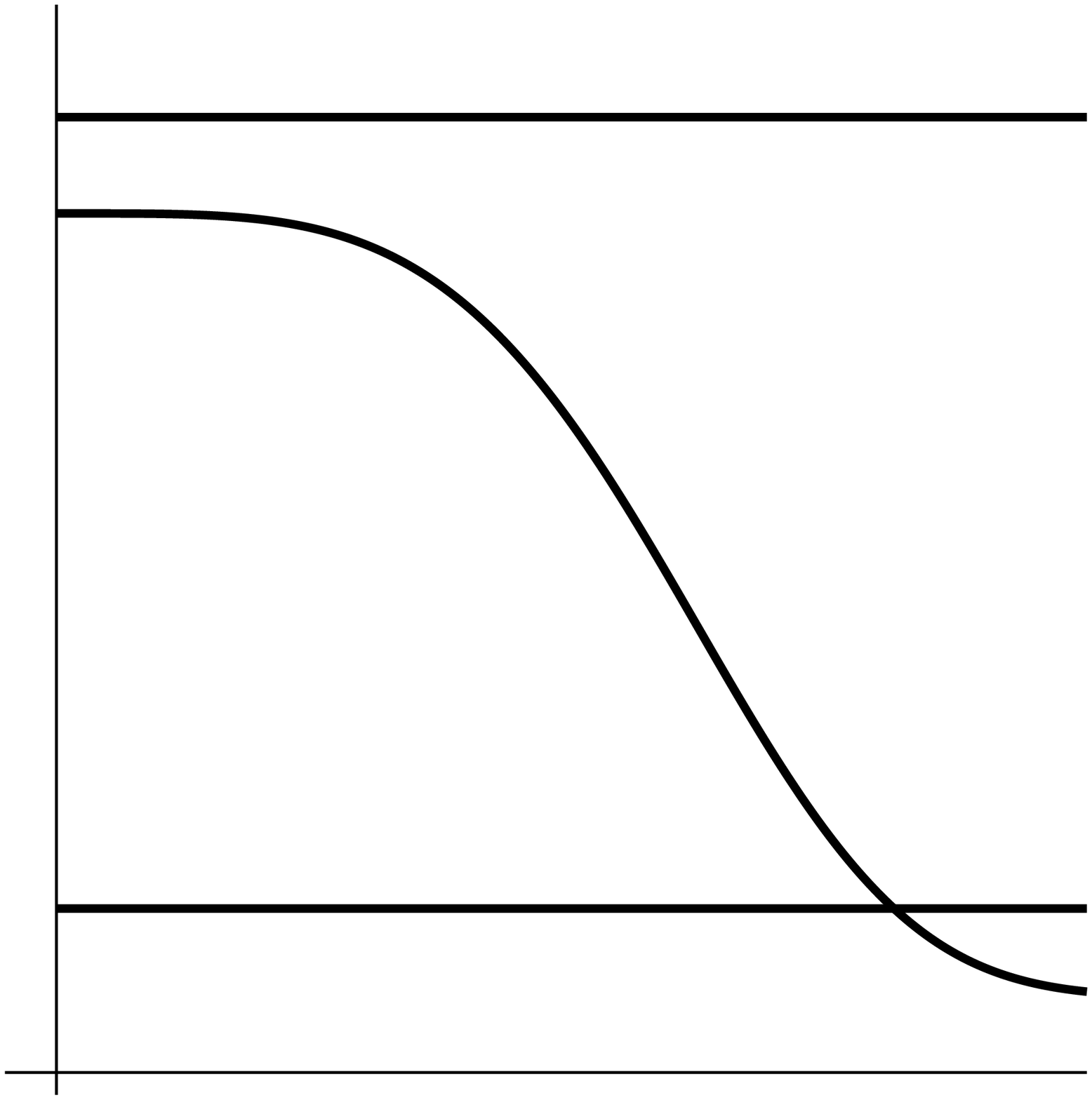}}
\put(0,25){\includegraphics[width=20mm]{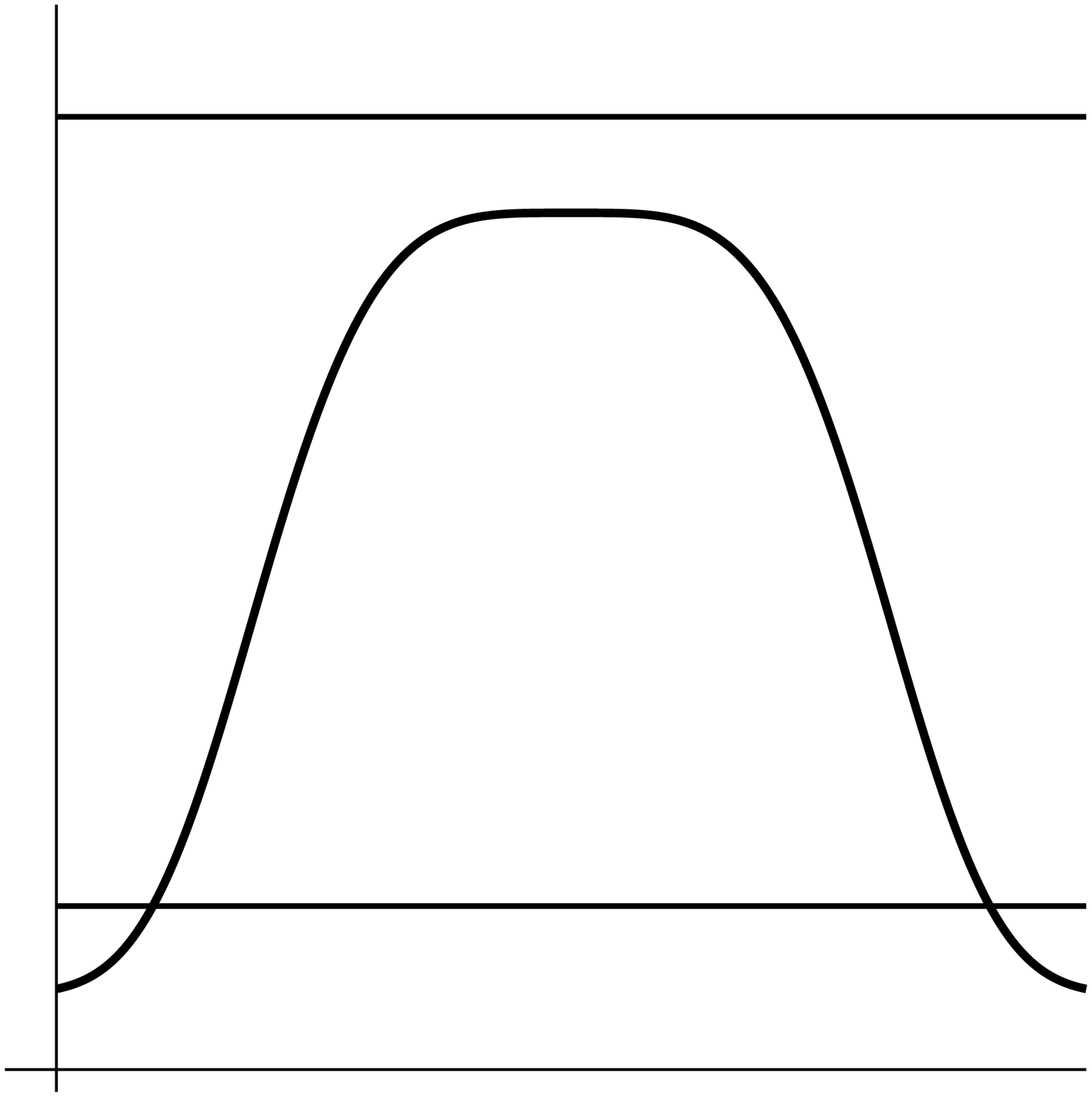}}
\put(26,25){\includegraphics[width=20mm]{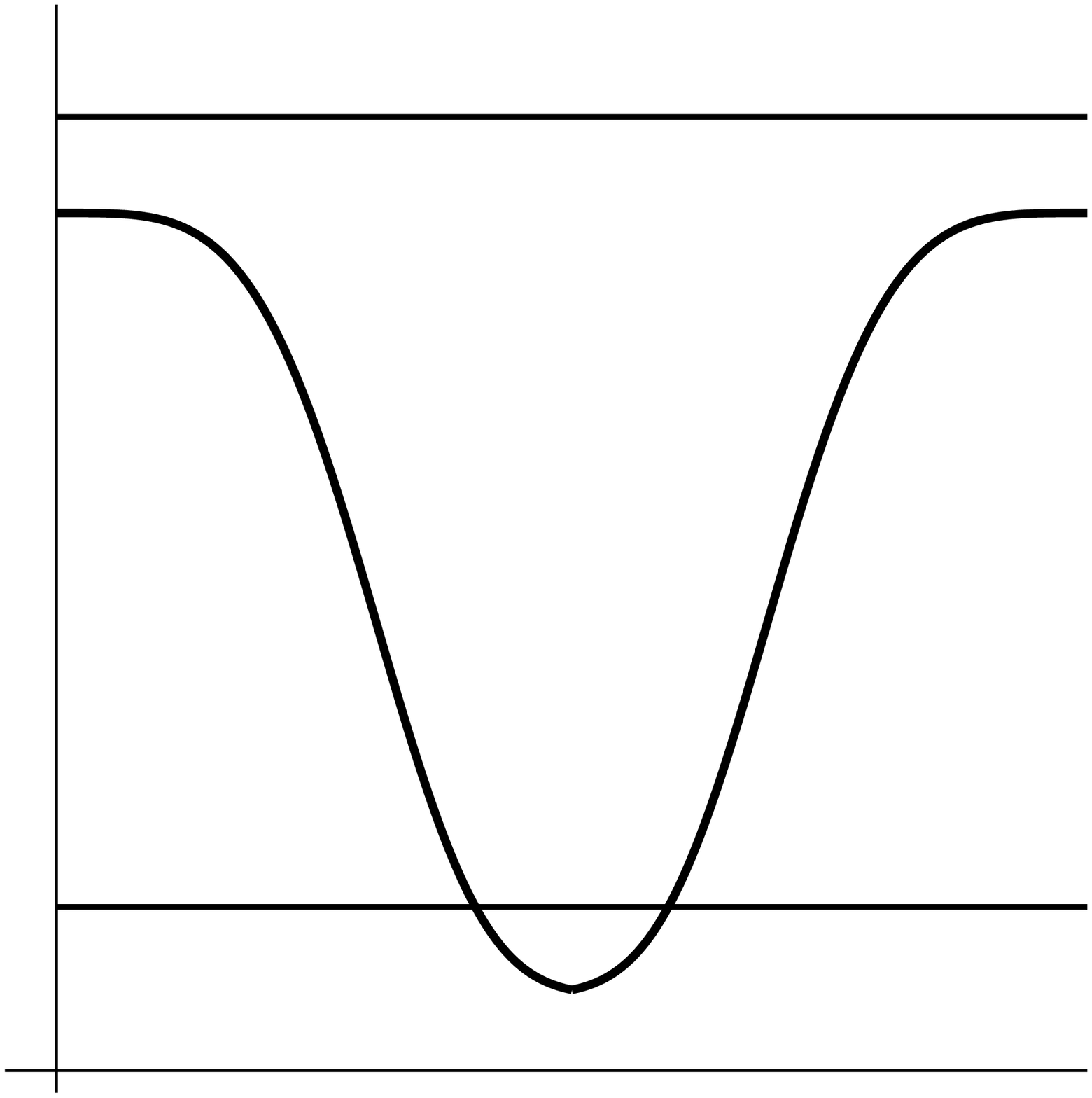}}
\put(0,0){\includegraphics[width=20mm]{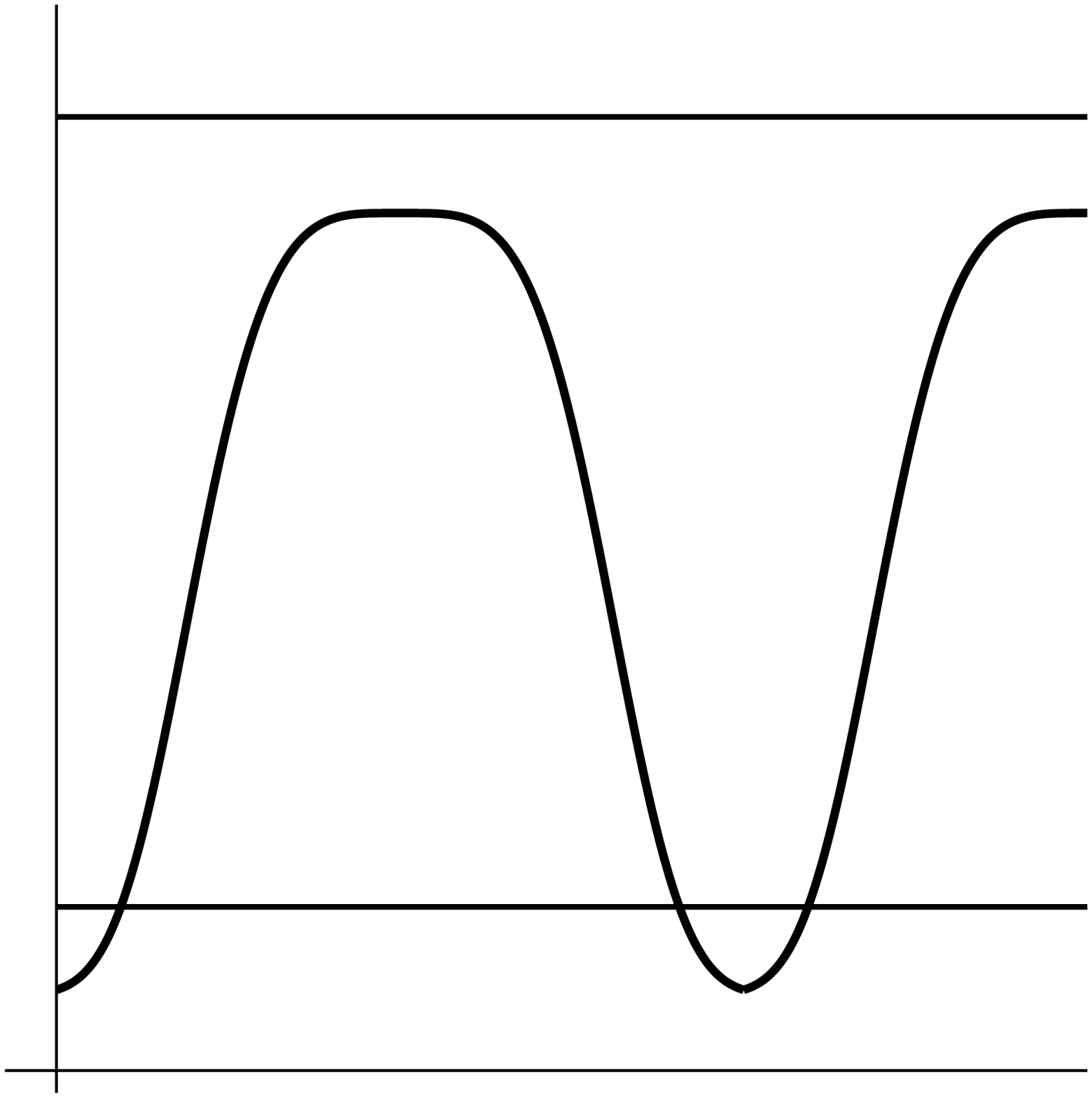}}
\put(26,0){\includegraphics[width=20mm]{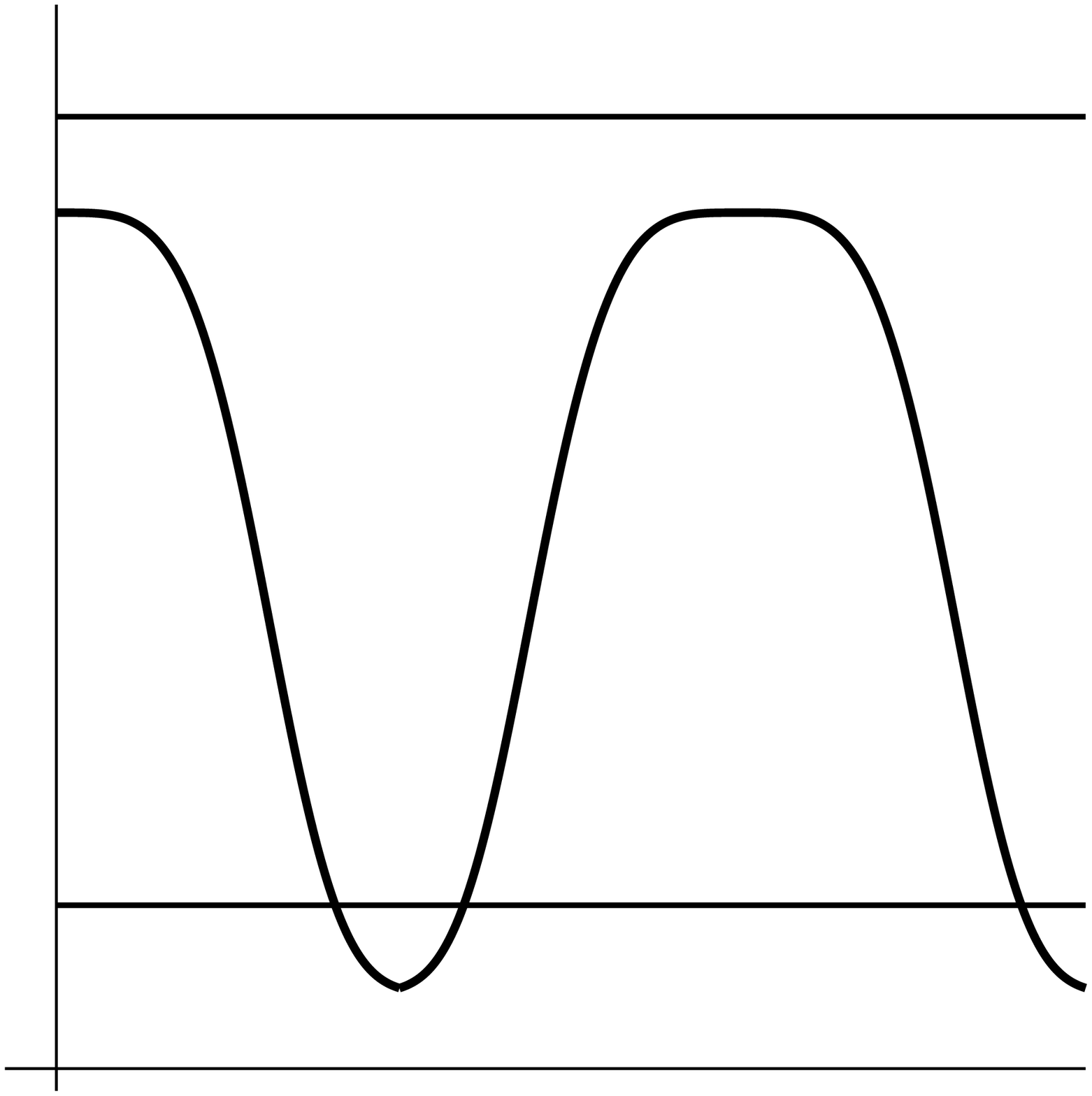}}

\put(-17,60){monotone}

\put(-15,34){2 modes}

\put(-15, 10){3 modes}

\put(-4, 17){$\bar w_{+}$}
\put(-4, 3){$\bar w_{-}$}

\put(-4, 42){$\bar w_{+}$}
\put(-4, 28){$\bar w_{-}$}

\put(-4, 67){$\bar w_{+}$}
\put(-4, 53){$\bar w_{-}$}

\put(22, 17){$\bar w_{+}$}
\put(22, 3){$\bar w_{-}$}

\put(22, 42){$\bar w_{+}$}
\put(22, 28){$\bar w_{-}$}

\put(22, 67){$\bar w_{+}$}
\put(22, 53){$\bar w_{-}$}

   \end{picture}
\caption{All positive solutions of problem \eqref{sw1}-\eqref{sw2} 
with $\gamma>0$ satisfying
 $3^2\gamma_0<\gamma<4^2\gamma_0$.}

  \label{fig:0}
}
\end{figure}


}

\begin{rem}
Using relations \eqref{s-uv},  we immediately obtain all solutions $(U(x), V(x), W(x))$ of system 
\eqref{seq1}-\eqref{sN}  such that $U(x)>0$ and $V(x)>0$ for all $x\in [0,1]$.
\end{rem}

\begin{rem}
It follows from the proof of Theorem \ref{thm:stationary} that 
non-constant 
solutions of problem \eqref{sw1}-\eqref{sw2}  oscillate around the constant solution
$\w_{-}$ and they are all below the constant $\w_{+}$.
\end{rem}

In our construction of non-constant solutions
 of the boundary value problem \eqref{sw1}-\eqref{sw2}, we use the well-known method from the 
classical mechanics that appears in the study of conservative systems with one degree of freedom (see {\it e.g.}
the Arnold book \cite{arnold}).
In Section  \ref{sec:pattern}, we recall that approach for reader's convenience. 
Our main contribution to this theory consists in the fact that we have found the optimal value of the coefficient $\gamma$ for which we 
have non-constant stationary solutions. Moreover, we prove that there is no other solutions, which is a consequence of the monotonicity of the function
$T=T(E)$ defined in \eqref{TE}, below.

Our next goal is to discuss nonnegative solutions  of \eqref{seq1}-\eqref{sN} such that $U(x)=V(x)=0$ 
on a certain set  $\I\subset [0,1]$ which we refer to  as a {\it null set} of a solution $(U,V,W)$.
 Here, we do not expect $W$ to be a $C^2$-function and, in the following, we say that a vector $(U,V,W)\in 
L^\infty ([0,1])\times L^\infty ([0,1]) \times C^1([0,1])$ is a {\it weak solution} of system \eqref{seq1}-\eqref{sN}
if the algebraic equations \eqref{seq1}-\eqref{seq2} are satisfied for almost all $x\in [0,1]$
and if
\begin{equation}\label{W:weak}
-\frac{1}{\gamma} \int_0^1 W'(x)\varphi'(x)\,dx
+   \int_0^1 \left( -d_g W(x) -U^2(x) W(x) +d V(x) +\kappa_0\right)\varphi(x)\,dx =0 
\end{equation}
for all $\varphi\in C^1([0,1])$.

Comparing to Theorem \ref{thm:stationary}, 
the  set of  weak  solutions is more complicated.

\begin{theorem}[Existence of discontinuous patterns] \label{thm:disc-stationary}
Assume that  $a>d_c$ and  $\kappa_0^2>\Theta$.  There exists a continuum of weak solutions 
of system \eqref{seq1}-\eqref{sN} with some $\gamma>0$. 
Each such a solution $(U,V,W)\in L^\infty(0,1)\times L^\infty(0,1) \times C^1([0,1])$
has the following property: 
there exists a sequence $0=x_0<x_1<x_2<...<x_N=1$ such that for each $k\in \{0, N-1\}$
either
\begin{itemize}
\item   for all $x\in (x_k,x_{k+1})$,
 $U(x)=V(x)=0$  and  $W(x)$ satisfies
$
\frac{1}{\gamma}W''-d_g W+\kappa_0=0 ,
$
\end{itemize}
or
\begin{itemize}
\item  for all $x\in (x_k,x_{k+1})$, $U(x)>0$ and $V(x)>0$  are given by relations \eqref{s-uv}, where
the function $W$ is a solution of  equation \eqref{sw1}.   
\end{itemize}
\end{theorem}

In Theorem \ref{thm:disc-stationary}, we do not attempt to classify all  discontinuous stationary solutions (as it was done in the continuous case in Theorem \ref{thm:stationary}), because they  appear to be  unstable
 solutions of the reaction diffusion equations \eqref{eq1}-\eqref{N}, see the next subsection.
Instead, in the proof of Theorem \ref{thm:disc-stationary}, we present a simple geometric argument which allows us to construct
 such discontinuous patterns.
Analogous constructions, using either geometrical or analytical arguments can be found in \cite{ATW88, MTH80}.

 \subsection{Instability of patterns.} \label{subsec:2.5}

Finally, we discuss  stability (in the sense of Lyapunov)
of stationary solutions of system \eqref{eq1}-\eqref{N}, which are constructed in Theorems \ref{thm:stationary}
and \ref{thm:disc-stationary}, 
and
surprisingly, we prove that they are all unstable.
Let us be more precise.

First, we consider a stationary solution  $(U(x), V(x), W(x))$ of  \eqref{eq1}-\eqref{N} where $W(x)$ is one of the
functions from Theorem \ref{thm:stationary} and $U(x)>0$, $V(x)>0$ are obtained from $W(x)$ via relations \eqref{s-uv}.
As the usual practice, we linearize system \eqref{eq1}-\eqref{eq3} around the steady state $(U(x), V(x), W(x))$
to obtain a system with three linear evolution equations  with the 
linear operator $\L$ defined formally by 
\begin{equation}\label{ELL0}
\L
\left(
\begin{array}{c}
 \varphi    \\
\psi   \\
\eta         
\end{array}
\right)
=
\left(
\begin{array}{ccc}
 0   &  0  &   0  \\
0     &       0       &  0\\
0         &    0     &   \frac{1}{\gamma} \partial_x^2\eta
\end{array}
\right)+\A
\left(
\begin{array}{c}
 \varphi    \\
\psi   \\
\eta         
\end{array}
\right),
\end{equation}
where 
\begin{equation}\label{A0}
\A(x) =(a_{ij}(x))_{i,j=1,2,3}\equiv 
\left(
\begin{array}{ccc}
 {d_c\left(\frac{d_c}{a}-1\right)}   &   \frac{(a-d_c)^2}{a}  &   0  \\
2 K                                                &         -d_b-d                     &   \frac{K^2}{W^2(x)}\\
-2 K                                              &              d                          &    -d_g- \frac{K^2}{W^2(x)}
\end{array}
\right),
\end{equation}
with the constant  
$
K\equiv U(x)W(x)= \frac{d_c(d_b+d)}{a-d_c}
$, see the beginning of Section \ref{sect:spect} for more details.
We consider $\L$ as a linear  operator in the Hilbert space 
$
\H = L^2(0,1)\times L^2(0,1)\times L^2(0,1)
$
with the domain
$
D(\L)= L^2(0,1)\times L^2(0,1)\times W^{2,2}(0,1).
$

It follows from direct calculations that the constant coefficient matrix
\begin{equation}\label{A12-i}
\A_{12}  \equiv 
\left(
\begin{array}{cc}
 {d_c\left(\frac{d_c}{a}-1\right)}   &   \frac{(a-d_c)^2}{a}    \\
2 K                                                &         -d_b-d                     
\end{array}
\right)
\end{equation}
obtained from $\A(x)$ after removing the last row and the last column,
has two real eigenvalues of opposite sign. The positive eigenvalue $\lambda_0$ 
of  $\A_{12}$  (given explicitly by formula \eqref{lambda_0})
plays the crucial role in the proof 
of instability of  stationary solution to problem \eqref{eq1}-\eqref{N}.

 \begin{theorem}[Instability of continuous patterns]\label{thm:unstab}
Consider the linear operator  $\L$ defined in \eqref{ELL0}, where
 $W(x)$ be one of the
functions from Theorem~\ref{thm:stationary} and $U(x)$, $V(x)$ are obtained from $W(x)$ via relations \eqref{s-uv},
except  the constant solution $(\u_+,\v_+,\w_+)$, see \eqref{const-pm}.
Then, the positive eigenvalue $\lambda_0$ of the matrix $\A_{12}$
belongs to the continuous spectrum of the operator $\L$.
Moreover, there exists a sequence $\{\lambda_n\}_{n\in \N}$
of positive eigenvalues of the operator $\L$ that satisfy $\lambda_n\to \lambda_0$ as $n\to \infty$.
 \end{theorem}

\begin{rem}
For simplicity of the exposition,
the constant stationary solution $(\u_+,\v_+,\w_+)$ is excluded from Theorem ~\ref{thm:unstab}, however, its 
instability is clear. Namely, it is an 
unstable solution of the kinetic system \eqref{heq1}--\eqref{heq3}, see Corollary  \ref{app-cor1}.
\end{rem}

Notice that, in order to present our idea in the simplest context, first, we prove Theorem~\ref{thm:unstab}
in the  case  of the constant stationary solution $(U(x),V(x),W(x))= (\u_-, \v_-, \w_-)$, see Theorem \ref{thm:instab:momog}, 
below. 

\begin{cor} \label{cor:unstab}
Every continuous stationary solution 
$(U, V, W)$  considered  in Theorem \ref{thm:unstab}
is an unstable solution of the nonlinear system \eqref{eq1}-\eqref{N}.
\end{cor}

This corollary   results immediately from the classical theory.
Indeed, by  Theorem \ref{thm:unstab}, every continuous stationary solution of system \eqref{eq1}-\eqref{N} 
is linearly unstable because the corresponding operator $\L$ has infinitely many positive eigenvalues.
Now, it suffices to note that $-\L$ is a sectorial operator, see the monograph by Henry \cite{Henry},
hence, applying the general result from \cite[Thm. 5.1.3]{Henry}  we show that $(U, V, W)$ 
is an unstable solution of  \eqref{eq1}-\eqref{N}.

Next, we show the instability of discontinuous stationary solutions from Theorem  \ref{thm:disc-stationary}. 
First, we have to emphasize that a direct linearization
of system \eqref{eq1}-\eqref{N} is not possible at a weak stationary solution satisfying
 $U(x)=V(x)=0$  at some $x\in [0,1]$. This is a consequence of the fact that the gradient of  the nonlinear term in equation
 \eqref{eq1} is not well-defined when $u=v=0$, see the Jacobian matrix \eqref{DF}, below.   
Hence, we  modify our  proof of instability in the following way.

 Let $(U_\I(x), V_\I(x), W_\I(x))$ be weak solution of  \eqref{seq1}-\eqref{sN} with a null  set $\I\subset [0,1]$,
namely, we assume that 
$$
 U_\I(x)=V_\I(x)=0 \;  \text{ for  } \;   x\in \I\quad   \text{and}   \quad  U_\I(x)>0 , \; V_\I(x)>0
\;  \text{ for  } \; x\in [0,1]\setminus \I.
$$  
Let us exclude the case $\I=[0,1] $, where the constant vector
$(U_\I, V_\I, W_\I)=(0,0, \kappa_0/d_g)$ is asymptotically stable solution of   \eqref{eq1}-\eqref{N} by Corollary \ref{cor:stab:trivial}.
For a null set $\I$, we define the associate $L^2$-space 
$$
L^2_\I(0,1)=\{ v\in L^2(0,1)\,:\, v(x)=0\quad \text{for}\quad x\in\I\}, 
$$
supplemented with the usual $L^2$-scalar product, which  is a Hilbert space as the closed subspace of $L^2(0,1)$.
Obviously,  when the measure of $\I$ equals zero, we have $L_\I^2(0,1)=L^2(0,1)$, thus, we assume in the following 
that $\I$ has
a positive Lebesgue measure and is different from the whole interval.

Now, observe that if $u_0(x)=v_0(x)=0$ for some $x\in [0,1]$ then by equations \eqref{eq1}-\eqref{eq2},
 we have $u(x,t)=v(x,t)=0$ for all $t\geq 0$. 
 Hence, in view of Theorem \ref{thm:existence2},  the space 
\begin{equation}\label{HI}
\H_\I=L^2_\I(0,1)\times L^2_\I(0,1) \times L^2(0,1)
\end{equation}
is invariant 
 for the flow generated by system \eqref{eq1}-\eqref{N} (notice that there is no ``$\I$'' in the last coordinate of $\H_\I$). 
 The crucial part of our analysis is based on the fact that, as long as we work in the space $\H_\I$,
 we can linearize system \eqref{eq1}-\eqref{N} at a discontinuous
 steady state $(U_\I, V_\I, W_\I)$. 
Moreover, for each $x\in [0,1]\setminus\I$, the corresponding linearized operator agrees with $\L$ defined in \eqref{ELL0}
 with the matrix $\A(x)$  from \eqref{A0}.
 Hence, the analysis from the proof of Theorem~\ref{thm:unstab} can be directly adapted to discontinuous steady states
 and we obtain the following counterpart of Corollary \ref{cor:unstab}.
 
\begin{cor} \label{cor:disc:unstab}
Every discontinuous weak stationary solution 
$(U_\I, V_\I, W_\I)$  of   problem \eqref{seq1}-\eqref{sN} with a null set $\I\subset [0,1]$,
which was constructed in Theorem \ref{thm:disc-stationary},
is an unstable solution of the nonlinear system \eqref{eq1}-\eqref{N},
considered in the Hilbert space $\H_\I$.
\end{cor}

\begin{rem}\label{rem:disc}
In other words, the  instability described by  Corollary \ref{cor:disc:unstab} appears when    
we perturb a sta\-tio\-na\-ry solution
 $(U_\I, V_\I, W_\I)$ on the set $[0,1]\setminus \I$, namely, in those points where $U_\I$ and $V_\I$ are nonzero.
On the other hand, by Theorem \ref{thm:stab:trivial},  small perturbations of
this steady state on  the null set $\I$, 
leads to a 
solution  satisfying  $u(x,t)\to 0$ and $v(x,t)\to 0$ as $t\to \infty$,  for each $x\in \I$.
\end{rem}

The proofs of Theorem  \ref{thm:unstab} 
and Corollary  \ref{cor:disc:unstab} are given in Section \ref{sect:spect}.
   
\subsection{Pattern formation}
A natural question arises if our model exhibits a formation of any  
pattern, which persist for long times. Numerical simulations indicating the growth of such, computationally stable, spatially heterogeneous solutions were performed in  \cite{MCK06, MCK07, MCK08}. The observed patterns take form of either periodic
or irregular spikes and,   now, it is clear that they cannot be described 
by stationary solutions, neither by smooth one from Theorem \ref{thm:stationary} nor
discontinuous one provided by Theorem \ref{thm:disc-stationary}.
Our present research is focused on understating these phenomena and answering questions on pattern formation in these kind of models. 


\section{Existence of solutions} \label{sec:existence}

We begin our study of  properties of solutions to the initial value problem \eqref{eq1}-\eqref{ini}
by showing  that it
  has  a unique and global-in-time solution
for all bounded and nonnegative initial conditions.

\begin{theorem}\label{thm:existence2}
Assume that  initial conditions 
 $u_0,v_0 \in L^\infty (0,1)$ and $w_0 \in W^{1, 2} (0, 1)$ are nonnegative. 
Then, for every $T>0$, the initial value problem 
\eqref{eq1}-\eqref{ini} has a unique, global-in-time, nonnegative solution 
$u \in C([0, T], L^\infty (0, 1))$, $v \in C([0, T], L^\infty (0, 1))$ and $w \in C ([0, T], W^{1, 2} (0, 1))$ such that %
$
u(x,\cdot),v(x,\cdot)\in C^1 ([0,T])
$  for every $x\in [0,1]$
and 
$
\;  w\in C^1 ( [0,T]; W^{1,2} (0,1)) \cap  C ( (0,T]; W^{2,2} (0,1)) .
$ 
\end{theorem}

\begin{rem}\label{thm:existence}
We can improve regularity of solutions using results from \cite[p.~112]{R84} in the following way. %
Assume that there exists $\alpha\in (0,1)$ such that 
 $u_0,v_0\in C^{\alpha} (0,1)$ for some $\alpha \in (0, 1)$ and if $w_0\in C^{2+\alpha} (0,1)$ satisfies $(w_0)_x(0)=(w_0)_x(1)=0$.
 Suppose, moreover, that 
 $u_0(x) > 0$, $v_0(x) > 0$, and $w_0(x) \ge 0$ for all $x\in [0,1]$.
Then, the initial value problem 
\eqref{eq1}-\eqref{ini} has a unique, global-in-time, and nonnegative solution 
$
u,v\in C^{\alpha,1+\alpha/2} ([0,1]\times [0,T]) \;  \text{and}
\;  w\in C^{2+\alpha,1+\alpha/2} ([0,1]\times [0,T])
$ 
for every $T>0$.
\end{rem}

The proof of Theorem \ref{thm:existence2} is more-or-less standard and we recall it for the completeness of the exposition.
 Below, we have only sketched the reasoning and we refer the reader to \cite[pp.~108-123]{R84} for other details.

First, we state a result on the local-in-time existence of solutions %
to the following initial-boundary value problem
\begin{align}
u_t &= f_1(u,v) \equiv \left(\frac{av}{|u| + |v|} - d_c \right)u, \label{local-eq1}\\
v_t &= f_2(u,v,w) \equiv -d_b v + u^2 w - dv, \label{local-eq2}\\
w_t &=\frac{1}{\gamma} w_{xx} -d_g w + f_3(u,v,w) \equiv \frac{1}{\gamma} w_{xx} -d_g w -u^2 w + dv + \kappa_0, \label{local-eq3}
\end{align}
supplemented with the homogeneous Neumann boundary conditions 
\begin{equation}\label{local-N}
w_x(0,t)=w_x(1,t)=0 \quad \mbox{for all} \quad t>0
\end{equation}
and with  initial conditions
\begin{equation}
u(x,0)=u_0(x), \quad v(x,0)=v_0(x), \quad w(x,0)=w_0(x).\label{local-ini}
\end{equation}

\begin{lemma}\label{lem:local}
Let $u_0, v_0 \in L^\infty (0, 1)$ and $w_0 \in W^{1, 2} (0, 1)$. There exists $T_0>0$ such that 
the initial value problem 
\eqref{local-eq1}-\eqref{local-ini} has a unique local-in-time solution 
$$
u \in C([0, T], L^\infty (0, 1)),\ v \in C([0, T], L^\infty (0, 1))\quad  \text{and}\quad  w \in C([0, T], W^{1, 2} (0, 1)).
$$ 
\end{lemma}

\begin{proof}
Here, it suffices to apply the abstract results by Rothe \cite[Thm.~1, p.~111]{R84}. %
We recall that a {\it mild} solution  $ (u,v,w)$ 
of problem \eqref{local-eq1}-\eqref{local-ini}, on a time interval 
$[0,T)$   and with initial data $u_0,v_0,w_0\in L^\infty (0,1)$,  
are  measurable functions $u,v,w :(0,1)\times (0,T) \to \R$ satisfying the following system of
integral equations  
\begin{align}
&u(x,t)= u_0(x)+ \int_0^t f_1(u(x,s),v(x,s))\,ds, \label{d-1}\\
&v(x,t)= v_0(x)+ \int_0^t f_2(u(x,s),v(x,s), w(x,s))\,ds,  \label{d-2} \\
&w(x,t)= S(t)w_0(x)+ \int_0^t S(t-s) (f_3(u,v,w))(x,s) \,ds, \label{d-3}
\end{align}
where $S(t)$ is a semigroup of linear operators associated with the equation $Z_t=\gamma^{-1} Z_{xx} - d_g Z$ %
on the interval $(0,1)$,  with the Neumann boundary conditions. %
Observe that the function $f_1 = f_1 (u, v)$ in equation \eqref{local-eq1} is Lipshitz continuous if we put $f_1 (0, 0) = 0$. %
Hence, to construct local-in-time solutions of system \eqref{d-1}-\eqref{d-3}, it suffices to apply the contraction mapping principle 
based on the Picard iterations. Details of such reasoning can be found in \cite[Thm.~1, p.~111]{R84}.
\end{proof}

Next, we show that the solutions are nonnegative.

\begin{lemma}\label{lem:positive}
Let the assumptions of Lemma \ref{lem:local} hold true.
Assume that $(u,v,w)$ is a solution of problem \eqref{local-eq1}-\eqref{local-ini}
on a certain interval $[0,T]$. %
If $u_0$, $v_0$ and $w_0$ are nonnegative, then $u(x, t)$, $v(x, t)$ and $w(x, t)$ are nonnegative %
for all $x\in [0,1]$ and $t\in [0,T]$.
\end{lemma}

\begin{proof}
{\it Step 1.}
First, notice that the inequality $w(x,t)\geq 0$ for all $x\in [0,1]$ and $t\in [0,T]$ implies 
$v(x,t) \ge 0$ for all $x\in [0,1]$ and $t\in [0,T]$. Indeed, using equation \eqref{local-eq2} with nonnegative $w$
we obtain the inequality $v_t(x,t)\geq -(d_b+d)v(x,t)$, which implies
$v(x,t)\geq e^{-(d_b+d)T}v_0(x)$ for all $x\in [0,1]$ and $t\in [0,T]$. %

{\it Step 2.}
Without loss of generality, we can assume that $w_0 (x) > 0$. Indeed, it suffices to replace $w_0$ by $w_0 + \varepsilon$ with %
$\varepsilon > 0$ and to use the continuous dependence of solutions of problem \eqref{local-eq1}--\eqref{local-ini} on inital conditions. %

For $w_0(x)>0$, by the continuity of $w$ (see Lemma \ref{lem:local}), we obtain $w(x,t)>0$ for all $x\in [0,1]$ and for sufficiently small $t>0$.
Suppose  there are $x_0\in [0,1]$ and  $t_0>0$ such that
$w(x,t)>0$ for $x\in [0,1]$ and $t\in [0,t_0)$ and such that $w(x_0,t_0)=0$. %
Note that $\sup_{0 \le 0 \le t_0}\|u (t)\|_\infty \le M$ for some constant $M > 0$ by Lemma \ref{lem:local} and %
$v(x, t) \ge 0$ for all $x \in [0, 1]$ and $t \in [0, t_0]$ by Step 1. %
Hence, using equation \eqref{local-eq3} we obtain 
\[
w_t \ge \frac{1}{\gamma} w_{xx} - d_g w - M^2 w + \kappa_0.
\]
In this inequality, all derivatives make sense, because, by a standard regularity theory for %
parabolic equation, we obtain
\[
w \in C^1 ([0, t_0], W^{1, 2}(0, 1)) \cap C ((0, t_0], W^{2, 2} (0, 1)).
\]
Since $\kappa_0 > 0$, by a comparison principle, we obtain $w(x, t) > 0$ for all $x \in [0, 1]$ and $t \in [0, t_0]$, %
which contradicts our hypothesis $w (x_0, t_0) = 0$.

{\it Step.~3.}
Let us finally find the lower bound for $u$. %
As long as $v = v(x, t)$ is nonnegative, it follows from equation \eqref{local-eq1} 
that $u_t(x,t)\geq -d_c u(x,t)$. Integrating this differential inequality we arrive at the estimate 
$u(x,t)\geq e^{-d_ct} u_0(x)$ for all $x\in [0,1]$ and $t\in [0,t_1)$.
\end{proof}

\begin{proof}[Proof of Theorem \ref{thm:existence2}]
Obviously, as long as a solution $(u, v, w)$ of \eqref{local-eq1}--\eqref{local-ini} is nonnegative, it is also a solution of our %
original problem \eqref{eq1}--\eqref{ini}. 

As the usual practice, to prove that local-in-time solutions from Lemma \ref{lem:local} exists for all $t>0$, it suffices to show
that for every $T>0$ 
$$
\sup_{t\in [0,T)} \big( \|u(t)\|_\infty+  
 \|v(t)\|_\infty+
 \|w(t)\|_\infty\big)<\infty,
$$
see {\it eg.} \cite[Thm.~1.iii, p.~111]{R84} for more details.

For positive $u$ and $v$, we have $v/(u+v)\leq 1$, hence, it follows from equation \eqref{eq1} that 
$u_t(x,t)\leq (a-d_c) u(x,t)$  for all $x\in [0,1]$ and $t\in [0, T)$. Integrating this differential inequality
we obtain
\begin{equation}\label{u:est:T}
\sup_{t\in [0,T)} \|u(t)\|_\infty \leq e^{(a-d_c)T}\|u_0\|_\infty \qquad \text{for every}\quad T>0.
\end{equation}

Next, we use the integral equation \eqref{d-2} for $v$. Computing the $L^\infty$-norm and using \eqref{u:est:T} we obtain
\begin{equation} \label{v:est:T}
\|v(t)\|_\infty \leq \|v_0\|_\infty +\int_0^t  \big((d_b+d)\|v(s)\|_\infty +C(T)\|w(s)\|_\infty\big)\,ds, 
\end{equation} 
where, following \eqref{u:est:T}, we denote $C(T)= e^{2(a-d_c)T}\|u_0\|_\infty^2$.

In a similar way, computing the $L^\infty$-norm of the integral equation \eqref{d-3} and using the well-known properties of %
the semigroup $S(t)$ we obtain
 \begin{equation} \label{w:est:T}
\|w(t)\|_\infty \leq \|w_0\|_\infty+\kappa_0 +\int_0^t  \big(d\|v(s)\|_\infty + C(T) \|w(s)\|_\infty\big)\,ds. 
\end{equation} 
Finally, adding inequalities \eqref{v:est:T} and \eqref{w:est:T} and using the Gronwall lemma, we complete the proof
that 
$
\sup_{t\in [0,T)} \big(  
 \|v(t)\|_\infty+
 \|w(t)\|_\infty\big)<\infty
 $
 for every $T>0$.
\end{proof}

\begin{proof}[Proof of Theorem \ref{thm:inv}.]
Integrating both sides of  equations \eqref{eq2} and \eqref{eq3} with respect to $x
\in [0, 1]$, adding the resulting formulas, and putting $\mu = \min
 \{d_g, d_b\}$, we obtain the differential inequality
\[
 \frac{d}{dt}\int_0^1 \big(v(x, t) + w(x, t)\big)\, dx \le -\mu \int_0^1 \big(v(x,
 t) + w(x, t)\big)\, dx + \kappa_0.
\]
Therefore, 
\begin{align}
 \int_0^1 \big(v(x, t) + w(x, t)\big)\, dx \le
 \left(\int_0^1 \big(v_0 (x) + w_0 (x)\big)\, dx \right) e^{-\mu t} +
 \frac{\kappa_0}{\mu}(1-e^{-\mu t}), \label{thm3.1-eq1}
\end{align}
which implies \eqref{limsup:2} and the bound by $\kappa_0/\mu$ in \eqref{limsup:1}, because $v, w$  are nonnegative.

To show the second bound in \eqref{limsup:1}, notice that, for nonnegative $u,v$, we have $u/(u+v)\leq 1$, hence, it follows from equation \eqref{eq1}
that $u_t\leq -d_c u+av$. Integrating this differential inequality, we obtain that the functions ${\U}(t)\equiv   \int_0^1 u(x,t)\;dx$ and 
${\V}(t)\equiv   \int_0^1 v(x,t)\;dx$
satisfy 
\begin{equation}\label{U:cal}
\U(t)\leq e^{-d_c t}\U(0)+a\int_0^t e^{-d_c(t-s)}\V(s)\;ds.
\end{equation}
Since $\V$ is a bounded function and $d_c>0$, we immediately obtain
$$
\int_0^{t/2} e^{-d_c(t-s)}\V(s)\;ds \leq e^{-d_ct/2}\int_0^{t/2} \V(s)\;ds \to 0\quad \text{as}\quad t\to\infty.
$$
On the other hand,  we have
$$
\int_{t/2}^t  e^{-d_c(t-s)}\V(s)\;ds \leq \left(\sup_{s\in [t/2,t]}\V(s)\right)
 \frac{1}{d_c} \left(1-e^{-d_ct/2}\right),
$$
where $\lim_{t\to\infty} \sup_{s\in [t/2,t]}\V(s) \leq \kappa_0/\mu$ by
 \eqref{limsup:2}. %
Hence, computing the limit superior as $t\to\infty$ 
of both sides of inequality \eqref{U:cal}, we complete the proof of ~\eqref{limsup:1}.

Next, let $B(x, y, t)$ be the  fundamental solution of the equation
$
Z_t = {\gamma}^{-1}Z_{xx} - d_g Z
$ 
 for $x \in (0, 1)$ and  $t >0$, supplemented with the Neumann boundary condition.
 It is well known (see {\it e.g. }\cite[pp.~19 \& 25] {R84})
 that $B(x, y, t)$ is nonnegative,  
\begin{align}
 \int_0^1 B(x, y, t)\, dy \le e^{-d_g t} \quad \text{for all}\ x \in [0,
 1], t > 0,\label{thm3.1-eq3.0}
\end{align}
and there exists a constant $C>0$ such that 
\begin{align}
 \int_0^1 B(x, y, t)v(y)\, dy \le Ct^{-1/2}e^{-d_g t} \|v\|_1 \label{thm3.1-eq3}
 \end{align}
 for each $v\in L^1(0,1)$, all  $x \in [0, 1] $, and  $t > 0$. %
Hence, equation \eqref{eq3} can be written as the following integral equation
\begin{align}
w(x, t) =& \int_0^1 B(x, y, t)w_0 (y)\, dy \label{w:integral} \\
&+ \int_0^t  \int_0^1 B(x, y,  t-s) \big( - u^2 (y, s)w(y, s) + d \,v(y,
 s)+\kappa_0 \big)   \, dy\, ds \nonumber
\end{align}
which, in view of the positivity of $B$, $w$, $v$,  via inequalities  \eqref{thm3.1-eq3.0},  implies
\begin{equation}
\begin{split}
\|w( t)\|_\infty \le \|w_0 \|_\infty e^{-d_g t} 
+ d \int_0^t  \int_0^1 B(x, y,
 t-s)v(y, s)\, dy \,ds + \frac{\kappa_0}{d_g}(1- e^{-d_g t}). 
\end{split}
\label{thm3.1-eq2}
\end{equation}
Moreover, it follows from \eqref{thm3.1-eq1} and \eqref{thm3.1-eq3} that the second term on the right-hand side of
\eqref{thm3.1-eq2} can be estimated as
\begin{equation}
\begin{split}
 d \int_0^t \int_0^1   B(x, y,
 t-s)&v(y, s)\, dy\,ds \\
\le&  
d C \int_0^t  (t-s)^{-1/2} e^{-d_g(t-s)}\|v( s)\|_1\, ds\\
\le &d C (\|v_0\|_1+\|w_0\|_1) \int_0^t (t-s)^{-1/2} e^{-d_g (t-s)}
 e^{-\mu s}\, ds \\
&+ dC \frac{\kappa_0}{\mu}\int_0^t  (t-s)^{-1/2} e^{-d_g (t-s)} \left(1-e^{-\mu
 s}\right)\, ds. \label{thm3.1-eq4}
\end{split}
\end{equation}
It is easy to prove that 
$
\lim_{t \to \infty} \int_0^t (t-s)^{-1/2}e^{-d_g (t-s)}
 e^{-\mu s}\, ds = 0
 $
 and
 $$
 \limsup_{t \to \infty}
 \int_0^t  (t-s)^{-1/2} e^{-d_g (t-s)} \left(1-e^{-\mu
 s}\right)\, ds = d_g^{-1/2} \int_0^\infty s^{-1/2} e^{-s}\,ds.
$$
Hence, the limit relation in  \eqref{limsup:3} is immediately obtained 
from \eqref{thm3.1-eq2} and \eqref{thm3.1-eq4}.
\end{proof}

\begin{rem}\label{bound-w-rem}
By an inspection of the proof of Theorem \ref{thm:inv}, we obtain the
 following global-in-time estimate of $w = w(x, t)$ :
\begin{align}
\|w(t)\|_\infty \le \|w_0\|_\infty + C d \left(2+\frac{1}{d_g}\right) \left(\|v_0 \|_1 +
 \|w_0 \|_1+ \frac{\kappa_0}{\mu}\right), \label{bound-w-eq1}
\end{align}
where $C$ and $\mu$ are the constants from \eqref{thm3.1-eq4}. %
\end{rem}

\section{Large time behavior of solutions } \label{sec:large}

In this section, we prove preliminary results on the large time behavior of solutions of problem \eqref{eq1}-\eqref{ini}.
Let us first prove that, under the assumption $a<d_c$, each positive solution of system \eqref{eq1}-\eqref{eq3} converges exponentially
towards the trivial steady state $(0,0, \kappa_0/d_g)$.

\begin{prop} \label{prop:a<dc}
Let $a < d_c$. Assume that $(u, v, w)$ is a nonnegative 
global-in-time solution of problem  \eqref{eq1}-\eqref{ini}  corresponding to
a bounded initial condition.
Then, 
there exist 
positive constants $C_1$, $C_2$, $\sigma_1$, $\sigma_2$  dependent  on the
parameters in system \eqref{eq1}-\eqref{ini} and 
 $C_2$ dependent  also on  $\|u_0 \|_\infty$  and 
$\|v_0\|_\infty$, such that
\begin{align}
0\leq u(x, t)& \le u_0 (x) e^{-(d_c - a)t}, \label{prop3.1-eq1}\\
0\leq v (x, t) &\le v_0 (x) e^{-(d_b + d)t}+ C_1 u_0^2 (x) t e^{-\sigma_1 t}, \label{prop3.1-eq2}\\
\left\|w (\cdot, t) - \frac{\kappa_0}{d_g}\right\|_{L^\infty} &\le C_2 e^{-\sigma_2
 t}, \label{prop3.1-eq3}
\end{align}
for all $x \in [0, 1]$ and $t \ge 0$.
\end{prop}

\begin{proof}
Since $u$ and $v$ are nonnegative, we have $v/(u + v) \le 1$. Hence, we
easily deduce from equation \eqref{eq1}  the following inequality %
$ u_t (x, t) \le -(d_c - a)u (x, t)$, %
which implies estimate \eqref{prop3.1-eq1}.

Since Theorem \ref{thm:inv} yields $K_w\equiv \sup_{t>0} \|w(t)\|_\infty<\infty$, then
using equation \eqref{eq2}  and estimate \eqref{prop3.1-eq1}, we obtain
 the differential inequality
\begin{align}
 v_t (x, t) \le -(d_b + d)v(x, t) + K_w u_0^2 (x) e^{-2(d_c - a)t}, \label{prop3.1-eq6}
\end{align}
which implies \eqref{prop3.1-eq2}.

To show the exponential convergence of $w=w(x, t)$ towards
 $\kappa_0/d_g$ stated in \eqref{prop3.1-eq3}, it suffices to use the
 integral representation \eqref{w:integral} of $w$, inequalities \eqref{prop3.1-eq1} and \eqref{prop3.1-eq2}, and to follow
 estimates from the proof of Corollary \ref{cor:stab:trivial}, below. %

\end{proof}


\begin{prop}\label{prop:zero}
Assume that $\big(u(x,t),v(x,t),w(x,t)\big)$ is a nonnegative solution of problem \eqref{eq1}-\eqref{ini}.
Fix $x\in [0,1]$. The following two conditions are equivalent
$$
(i)  \; \;  u(x,t)\to 0  \quad \text{as}\quad  t\to\infty  \qquad \text{and}\qquad 
(ii)\; \;  v(x,t)\to 0   \quad \text{as} \quad  t\to\infty.
$$
\end{prop}

\begin{proof}
We rewrite equation \eqref{eq2} in the following integral form
$$
v(x,t) = e^{-(d_b+d)t}v_0(x) + \int_0^t e^{-(d_b+d)(t-s)}u^2(x,s)w(x,s)\;ds.
$$ 
By Theorem \ref{thm:inv}, the function $w(x,t)$ is bounded for $x\in [0,1]$ and $t\geq 0$, hence, combining 
condition (i) with the Lebesgue dominated convergence theorem in this integral equation, we can prove (ii).
 
Recall that for nonnegative $u$ and $v$ we have $u/(u+v)\leq 1$. Hence, it follows from equation \eqref{eq1} that
$u_t(x,t)\leq -d_c u(x,t) +av(x,t)$. Integrating this differential inequality and following the previous argument %
we prove directly that  (ii) implies  (i).
\end{proof}

Before proving Theorem \ref{thm:stab:trivial}, 
in the following lemma, we show a kind of pointwise stability of the trivial steady state.

\begin{lemma}\label{lem:stab}
Let $x\in [0,1]$.
Under the assumptions of Theorem \ref{thm:stab:trivial}, it holds 
\begin{equation}\label{lem:uv}
0\leq u(x,t) < M\left(1+\frac{d_c}{a}\right) \qquad \text{and}\qquad 0\leq v(x,t)<\left(\frac{d_c}{a}\right)^2 M\\ 
\end{equation}
for all $t>0$.
\end{lemma}


\begin{proof}
Recall that the nonnegativity of $u$ and $v$ has been shown already in Lemma  \ref{lem:positive}.

{\it Step 1.} First, we suppose that there is $T>0$ such that $u(x,t)\leq M(1+d_c/a)$ for all $t\in [0,T]$.
Remember that $0\leq w(x,t)\leq K_w$ for all $x\in [0,1]$ and $t\geq 0$. Hence, 
it follows from equation \eqref{eq2} satisfied by  $v$ that 
\begin{equation}\label{v:M:est}
\begin{split}
\frac{\partial}{\partial t}  \left( v(x,t)-\left(\frac{d_c}{a}\right)^2 M\right)= &
- (d_b+d) \left( v(x,t)-\left(\frac{d_c}{a}\right)^2 M\right) +u^2w\\
& -(d_b+d) \left(\frac{d_c}{a}\right)^2 M\\
\leq &
- (d_b+d) \left( v(x,t)-\left(\frac{d_c}{a}\right)^2 M\right) \\
& +\left[ M^2  \left(1+\frac{d_c}{a}\right)^2  K_w-(d_b+d) \left(\frac{d_c}{a}\right)^2 M\right].
\end{split}
\end{equation}
The term in the brackets on the right hand side of \eqref{v:M:est} is non-positive, because of  assumption \eqref{as:M},
hence, we skip it. Integrating the resulting differential inequality,  we obtain
\begin{equation}\label{v:M:est:fin}
v(x,t) -\left(\frac{d_c}{a}\right)^2 M \leq e^{-(d_b+d) t} \left( v_0(x) -\left(\frac{d_c}{a}\right)^2 M \right)
\end{equation}
for all $t\in [0,T]$. 
Using the assumption of $v_0(x)$ stated in \eqref{u:M}, we see that the right-hand side of \eqref{v:M:est:fin} is negative,
hence, we obtain 
$$
v(x,t) < \left(\frac{d_c}{a} \right)^2 M \qquad \text{for all} \quad t\in [0,T].
$$

{\it Step 2.} Now,   suppose that there is $T>0$ such that  $v(x,t) \leq (d_c/a)^2 M$ for all $t\in [0,T]$.
Applying 
  the inequality $u/(u+v)\leq 1$ in equation \eqref{eq1} yields
$$
u_t(x,t) \leq av(x,t) -d_c u(x,t) \leq a\left(\frac{d_c}{a}\right)^2 M-d_cu(x,t).
$$
Integrating this differential inequality and using the assumption on $u_0(x)$ from \eqref{u:M}, we obtain
\begin{equation}\label{lem:u:est}
u(x,t)\leq  e^{-d_c t} u_0(x) + \frac{a}{d_c} \left(\frac{d_c}{a}\right)^2 M \left(1-e^{-d_ct}\right)
< M\left(1+\frac{d_c}{a} \right) 
\end{equation}
for all  $t\in [0,T]$.

{\it Conclusion.} By Theorem \ref{thm:existence2},  for each $x\in [0,1]$, the functions
$u(x,t)$ and $v(x,t)$ are continuous with respect to  $t$. Suppose that inequalities
\eqref{lem:uv} hold true for all $t\in [0,T)$ with some $T>0$ and at least one of them becomes an equality for $t=T$. 
Such a hypothesis, however, contradicts either the implication from Step 1 or from Step 2.
This completes the proof of Lemma \ref{lem:stab}.
\end{proof}

\begin{proof}[Proof of Theorem \ref{thm:stab:trivial}]
We are going to construct 
 a sequence $\{T_{2n}\}_{n \ge 1}$, independent of $x\in [0,1]$, satisfying  
$T_{2n}<T_{2(n+1)}$  and $\lim_{n\to\infty} T_{2n}=+\infty$,
such that $u(x, t)$ and $v(x, t)$ fulfill   for each $n\in \N$ the following inequalities
\begin{align}
u(x, t) \le \theta^{n-1}\dfrac{d_c}{a}\left(1 + \dfrac{d_c}{a}\right)M
\quad \text{and}\quad v(x, t) \le
 \theta^{n}\left(\dfrac{d_c}{a}\right)^2 M \qquad \text{for all}\quad  t \ge
 T_{2n}. \label{th2.3-eq1}
\end{align}
Here $\theta \equiv [1+(d_c / a)^2]/2 < 1$, because of the assumption $d_c<a$. 
Since $u$ and $v$ are nonnegative, inequalities \eqref{th2.3-eq1}
imply
$
 \lim_{t \to \infty} u(x, t) = \lim_{t \to
\infty} v(x, t) = 0.
$

We proceed by induction with respect to $n$. %

{\it Step 1.}
By Lemma \ref{lem:stab}, we have $v(x, t) <
(d_c /a)^2 M$   for all $t>0$. The first inequality in \eqref{lem:u:est} and the assumption on $u_0(x)$ lead to the estimates
\begin{align}\label{th2.3-eq2} 
u(x, t) \le e^{-d_c t}u_0 (x) + \dfrac{a}{d_c}
 \left(\dfrac{d_c}{a}\right)^2 M (1-e^{-d_c t})
 \le e^{-d_c t}M + \frac{d_c}{a}M (1-e^{-d_c t}). 
\end{align}
Hence, choosing  $T_1 > 0$ such that 
$
 e^{-d_c t} M \le \left(\frac{d_c}{a}\right)^2 M 
$
for all  $t \ge T_1$
 we obtain
\begin{align}
u(x, t) \le \dfrac{d_c}{a} M \left(1 + \dfrac{d_c}{a}\right)\qquad
 \text{for all}\quad  t \ge T_1. \label{th2.3-eq3}
\end{align}
Consequently, it follows from the differential equation \eqref{eq2} and from the assumption on $w(x,t)$ 
that 
\begin{align}
v_t(x, t) \le -(d_b + d) v + \left[\dfrac{d_c}{a} M \left(1 +
 \dfrac{d_c}{a}\right)\right]^2 K_w \qquad \text{for all}\quad  t \ge T_1. \label{th2.3-eq4}
\end{align}
Integrating this differential inequality for  $t\geq T_1$, we
obtain
\begin{align}
v(x, t) \le e^{-(d_b + d)(t-T_1)}v(x, T_1) +
 \dfrac{\left(\frac{d_c}{a}\right)^2 M^2 \left(1 +
 \frac{d_c}{a}\right)^2 K_w}{d_b + d}\left(1 - e^{-(d_b +
 d)(t-T_1)}\right). \label{th2.3-eq5}
\end{align}
Using the assumption \eqref{as:M},  written in the form 
\[
 \dfrac{\left(\frac{d_c}{a}\right)^2 M^2 \left(1 +
 \frac{d_c}{a}\right)^2 K_w}{d_b + d} \le \left(\frac{d_c}{a}\right)^4 M,
\]
and the  inequality  $v(x, T_1) \le (d_c /a )^2
M$ from Lemma \ref{lem:stab},
 we obtain from \eqref{th2.3-eq5} that
\begin{align}
v(x, t) \le e^{-(d_b + d)(t-T_1)}\left(\dfrac{d_c}{a}\right)^2 M +
 \left(\dfrac{d_c}{a}\right)^4 M \qquad \text{for all}\quad  t \ge T_1. \label{th2.3-eq6}
\end{align}
Now, we choose  $T_2 > T_1$ such that 
\[
 e^{-(d_b + d)(t-T_1)}\left(\dfrac{d_c}{a}\right)^2 M \le
 \frac{1}{2}\left(\left(\dfrac{d_c}{a}\right)^2 M -
 \left(\dfrac{d_c}{a}\right)^4 M \right) \qquad \text{for all}\quad t \ge T_2.
\]
Therefore, the right-hand side of \eqref{th2.3-eq6} is estimated as 
\begin{align}
v(x, t) \le \frac{1}{2}\left(1-\left(\frac{d_c}{a}\right)^2
 \right)\left(\frac{d_c}{a}\right)^2 M + \left(\frac{d_c}{a}\right)^4 M
 \label{th2.3-eq7}
= \left(\frac{d_c}{a}\right)^2 M \theta \qquad \text{for all}\quad t \ge T_2, \nonumber
\end{align}
where $\theta \equiv  [1 + (d_c /a)^2]/2 < 1$. 
Consequently, we have found $T_2>0$ and we have proved the inequalities 
\[
 u(x, t) \le \frac{d_c}{a}M \left(1 + \frac{d_c}{a}\right) \quad \text{and}\quad  v(x,
 t) \le \left(\frac{d_c}{a}\right)^2 M \theta \qquad \text{for all}\quad t \ge
 T_2,
\]
which correspond to  \eqref{th2.3-eq1}  when $n = 1$.

{\it Step 2}. 
Now, we assume that \eqref{th2.3-eq1} holds true for $n\in \N$ and we prove it for $n+1$.
Equation~\eqref{eq1} together with  the inductive hypothesis for $v(x,t)$ lead to 
\begin{align}\label{th2.3-eq9}
u_t \le a v(x, t) - d_c u(x, t)
\le a \theta^n \left(\frac{d_c}{a}\right)^2 M - d_c u(x, t) \qquad
 \text{for all}\quad t \ge T_{2n}. 
\end{align}
Integrating the differential inequality  \eqref{th2.3-eq9} for $t\geq T_{2k}$, we
obtain 
\begin{equation} \label{th2.3-eq10}
\begin{split}
u (x, t) \le & e^{-d_c (t - T_{2n})} u(x, T_{2n}) +
 \dfrac{a}{d_c}\theta^n \left(\frac{d_c}{a}\right)^2 M \left(1-e^{-d_c
 (t - T_{2n})}\right) \\
\le& e^{-d_c (t - T_{2n})} \theta^{n-1}\frac{d_c}{a}\left(1 +
 \frac{d_c}{a}\right)M \\
&+
 \dfrac{a}{d_c}\theta^n \left(\frac{d_c}{a}\right)^2 M \left(1-e^{-d_c
 (t - T_{2n})}\right) \qquad \text{for all}\quad  t \ge T_{2n}. 
\end{split}
\end{equation}
Now, we choose  $T_{2n + 1} > T_{2n}$ in such a way  that that following inequality holds
\[
 e^{-d_c (t - T_{2n})} \theta^{n-1}\frac{d_c}{a}\left(1 +
 \frac{d_c}{a}\right)M \le \left(\frac{d_c}{a} \theta^n \right)^2 M
 \qquad \text{for all}\quad  t \ge T_{2n + 1},
\]
which  is equaivalent to the following one
$$ e^{-d_c (t - T_{2n})} (1 + (d_c /a)) \le (d_c /a)\theta^{n+1} 
\qquad \text{for all}\quad  t \ge T_{2n}.
$$ 
Hence, for $t \ge T_{2n + 1}$ we obtain from \eqref{th2.3-eq10} that 
\begin{align}\label{th2.3-eq11}
u(x, t) \le \left(\frac{d_c}{a} \theta^n \right)^2 M + \theta^n
 \dfrac{d_c}{a} M 
=  \dfrac{d_c}{a} M \theta^n \left(1 +  \dfrac{d_c}{a} \right)\dfrac{1
 +  \frac{d_c}{a} \theta^n}{1 +  \frac{d_c}{a}}.
\end{align}
Since $(1 + (d_c /a) \theta^n)/ (1 + (d_c /a)) < 1$, we have
\begin{align}
u(x, t) \le \dfrac{d_c}{a} M \theta^n \left(1 + \dfrac{d_c}{a}\right)
 \qquad \text{for all}\quad  t \ge T_{2n+1}. \label{th2.3-eq12}
\end{align}

Next,  integrating equation \eqref{eq1} and using \eqref{th2.3-eq12}, we obtain 
 that the function 	$v(x, t)$ satisfies for  $t \ge T_{2n+1}$ the following inequality
\begin{equation}\label{th2.3-eq13}
\begin{split}
v(x, t) \le& e^{-(d_b + d)(t-T_{2n+1})}v(x, T_{2n+1})\\& +
 \dfrac{\left(\frac{d_c}{a} M \theta^n \left(1 +
 \frac{d_c}{a}\right)\right)^2 K_w}{d_b + d} \left(1 - e^{-(d_b +
 d)(t-T_{2n+1})}\right). 
\end{split}
\end{equation}
It follows from the assumption \eqref{as:M} that 
\[
 \dfrac{\left[\frac{d_c}{a} M \theta^n \left(1 +
 \frac{d_c}{a}\right)\right]^2 K_w}{d_b + d} \le
 \left(\frac{d_c}{a}\right)^4 M \theta^{2n}
\]
and from the inductive hypothesis \eqref{th2.3-eq1} that $v(x, T_{2n+1}) \le \theta^n (d_c/a)^2 M$. %
Hence, choosing $T_{2(n+1)} > T_{2n + 1}$ in such a way that
\[
 e^{-(d_b + d)(t-T_{2n+1})} \theta^n \left(\frac{d_c}{a} \right)^2 M \le
 \frac{1}{2}\left[\left(\frac{d_c}{a} \right)^2 \theta^n  M -
 \left(\frac{d_c}{a} \right)^4 \theta^{2n} M \right] \quad \text{for all}\quad  t
 \ge T_{2(n+1)}
\]
or equivalently that  
\[
 e^{-(d_b + d)(t-T_{2n+1})} \le \frac{1}{2}\left(1 - \left(\frac{d_c}{a}
 \right)^2 \theta^n \right) \qquad \text{for all}\quad  t
 \ge T_{2(n+1)},
\]
we deduce  from \eqref{th2.3-eq13} that %
\begin{align}
v(x, t) \leq  \dfrac{1}{2}\left(\dfrac{d_c}{a} \right)^2 \theta^n M \left(1 +
 \left(\dfrac{d_c}{a} \right)^2 \theta^n \right) \qquad \text{for all}\quad  t
 \ge T_{2(n+1)}.\label{th2.3-eq13b}
\end{align}
Since $[1 + (d_c/a)^2 \theta^n]/ 2 < \theta$, inequality \eqref{th2.3-eq13b} implies 
\begin{align}
v(x, t) \le  \left(\dfrac{d_c}{a} \right)^2 M \theta^{n+1}\qquad
 \text{for all}\quad  t \ge T_{2(n+1)}. \label{th2.3-eq114}
\end{align}

Consequently,  inequalities  \eqref{th2.3-eq12} and \eqref{th2.3-eq114} hold true for all $t \ge T_{2(n+1)}$
and the proof of the inductive step is complete.
\end{proof}

\begin{proof}[Proof of Corollary \ref{cor:stab:trivial}.]
It follows immediately from Theorem  \ref{thm:stab:trivial} 
that $u(x,t)\to 0$ and $v(x,t)\to 0$ as $t\to\infty$ uniformly in $x\in [0,1]$.
To show that $\|w(t)-\kappa_0/d_g\|_\infty \to 0$ as $t\to \infty$ we
 use the fundamental solution $B(x, y, t)$ of the equation
$
Z_t = {\gamma}^{-1}Z_{xx} - d_g Z
$ 
for $x \in (0, 1)$ and  $t >0$, supplemented with the Neumann boundary
 condition. 
Noting \eqref{thm3.1-eq3.0} and the positivity of all functions
in the integral representation of $w$ from \eqref{w:integral} 
we obtain the
 inequality ({\it c.f.}  \eqref{thm3.1-eq2})
\begin{align}
\begin{split}
\|w(t) - \kappa_0/d_g \|_\infty \le \|w_0 - \kappa_0/d_g \|_\infty e^{-d_g t} 
+ d \int_0^t  \int_0^1 B(x, y,
 t-s)v(y, s)\, dy \,ds.
\end{split}\label{proof_cor2.4_eq1}
\end{align}
Since $v(x, t) \le (d_c/a)^2 M$ for all $x \in [0, 1]$ and $t > 0$ by Lemma \ref{lem:stab}, we see
\begin{align*}
\int_0^{t/2} \int_0^1
 B(x, y, t-s)v(y, s)\, dy \,ds \le \left(\frac{d_c}{a}\right)^2 M
 \int_0^{t/2} e^{-d_g (t-s)}\, ds \ \to 0\quad \text{as}\quad t\to\infty.
\end{align*}
On the other hand, 
\begin{align*}
\int_{t/2}^t \int_0^1 B(x, y, t-s)v(y, s)\, dy \,ds \le
 \left(\sup_{s\in [t/2,t]} v(x, s)\right)
 \frac{1}{d_g} \left(1-e^{-d_g t/2}\right).
\end{align*}
It follows from Theorem \ref{thm:stab:trivial} that the right-hand side
 above tends to $0$ as $t \to \infty$. Consequently, we conclude from
 \eqref{proof_cor2.4_eq1} that $\|w(t) - \kappa_0/d_g \|_\infty \to 0$
 as $t \to \infty$. 
\end{proof}

\section{Construction of patterns}\label{sec:pattern}

In this section, we discuss stationary solutions of system \eqref{eq1}-\eqref{N} (or equivalently of system
\eqref{seq1}-\eqref{sN}). 
First,  we limit ourselves to nonzero $U(x)$ and $V(x)$ to  obtain relations
\eqref{s-uv} as well as the boundary value problem \eqref{sw1}-\eqref{sw2}.

Let us begin with  
an arbitrary $C^1$-function $h:(0,\infty)\to\R$.
By the change of variables
\begin{equation} \label{change}
x\mapsto Tx, \qquad \text{where}\qquad T=\sqrt{\gamma},
\end{equation}
one can transform  the boundary value problem 
$$
\frac{1}{\gamma} w''+h(w)=0, \quad w'(0)=w'(1)=0, \qquad\text{with $x\in (0,1),$} 
$$
into the problem
\begin{eqnarray}
&& w''+h(w)=0, \qquad x\in (0,T), \label{BV1}\\
&&w'(0)=w'(T)=0. \label{BV2} 
\end{eqnarray}
Together with equation \eqref{BV1}, we consider the corresponding system of the first order equations 
\begin{eqnarray}
w'=z, \quad 
z'=-h(w).\label{sys2}
\end{eqnarray}
Hence,  solutions of equation \eqref{BV1} satisfying the boundary conditions \eqref{BV2}
correspond to  trajectories of system \eqref{sys2} satisfying $z(0)=z(T)$  for a certain $T>0$.
Here, we should recall that system  \eqref{sys2} in  autonomous in the following sense: if $(w(x),z(x))$
is  a solution, then  $(w(x+x_0),z(x+x_0))$ is a solution for all $x_0\in\R$, as well.

Multiplying  second equation in \eqref{sys2} by $z=z(x)$ and using first one, we obtain the equation for all trajectories 
of system \eqref{sys2}
\begin{equation}
\frac{z^2}{2} +H(w)=E, \label{energy}
\end{equation}
where $H'=h$ 
  and $E\in\R$ is an arbitrary  constant. Recall that  the constant $E$ is called  the {\it total energy} 
in the classical mechanics,  the function $H$ corresponds to the {\it potential energy}, and relation \eqref{energy} 
is the {\it first integral}  of system \eqref{sys2}.

 It follows from equation \eqref{energy} that all trajectories of system 
\eqref{sys2} are symmetric with respect to the $w$-axis. Hence, the condition $z(0)=w'(0)=z(T)=w'(T)$ 
is satisfied for a certain 
$T>0$ if equation \eqref{energy} describes a closed curve on the $wz$-plane for a some $E\in\R$. Such  closed curves exist only 
if the potential energy $H=H(w)$ 
has a local minimum at a certain point $\w_{-}$, see Fig. \ref{fig:2}.

%

\nofig{
%
%
%

\begin{figure}
  \setlength{\unitlength}{1mm}
{\footnotesize
\begin{picture}(60,90)
    \put(5,42){\includegraphics[width=60mm]{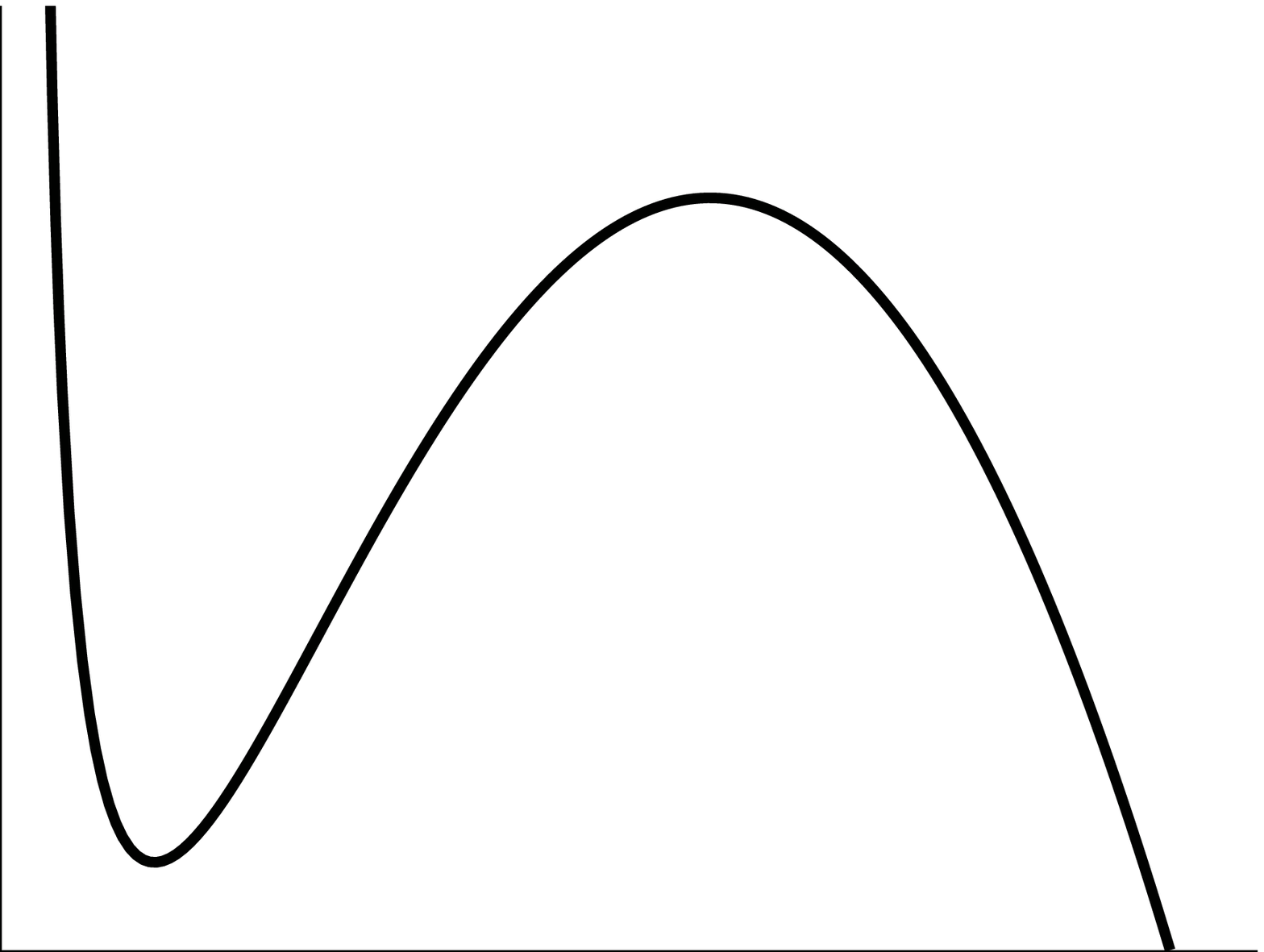}}
    \put(5,0){\includegraphics[width=60mm]{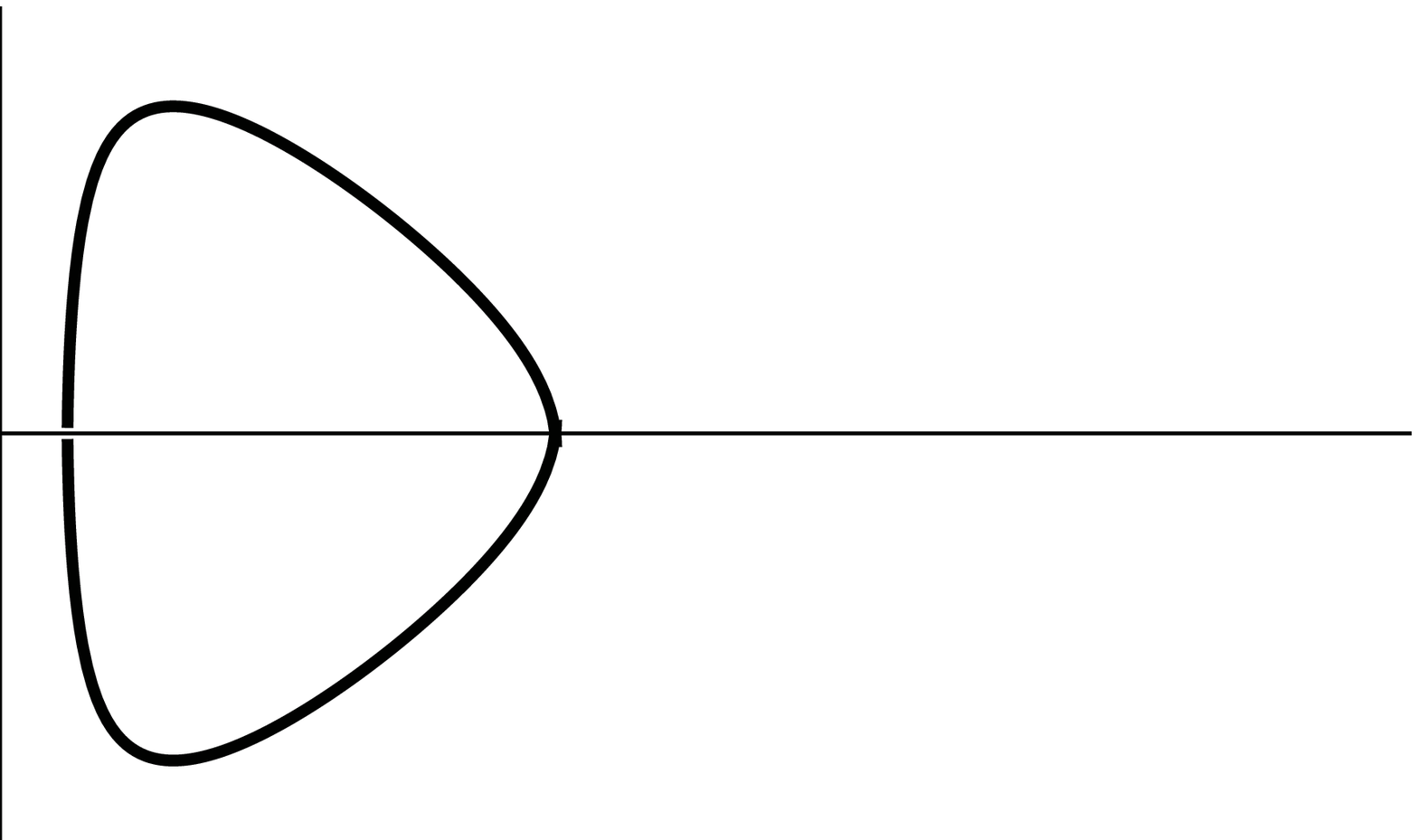}}
    \put(2,33){$Z$}
    \put(61,14){$W$}
    \put(2,69.5){$E$}
    \multiput(5,70.5)(2,0){12}{\line(1,0){1}}
    \put(-7,77){$H(\bar{w}_{+})$}
    \put(-7,45){$H(\bar{w}_{-})$}
    \put(36,39){$\bar{w}_{+}$}
    \put(12,39){$\bar{w}_{-}$}
    \put(3,39){$w_{1,E}$}
     \put(26,39){$w_{2,E}$}
   \multiput(5,78)(2,0){17}{\line(1,0){1}}
    \multiput(5,46)(2,0){4}{\line(1,0){1}}
    \put(58,54){$H=H(w)$}   
    \put(25,9){$z=\pm\sqrt{2(E-H(w))}$}  

    \multiput(7.5,20)(0,2){9}{\line(0,1){1}} 
    \multiput(28.5,20)(0,2){9}{\line(0,1){1}} 

    \multiput(7.5,41)(0,2){15}{\line(0,1){1}} 
    \multiput(28.5,41)(0,2){15}{\line(0,1){1}}

\multiput(38,42)(0,2){18}{\line(0,1){1}} 
\multiput(12.5,42)(0,2){2}{\line(0,1){1}} 
  \end{picture}
  \caption{The graph of the potential energy $H=H(w)$ defined in 
\eqref{H} and the closed trajectory of system \eqref{sys2} with the energy level $E>0$.}
  \label{fig:2}
}
\end{figure}

%
%
}


Let us  apply these   classical ideas to obtain preliminary results on the nonexistence of stationary 
solutions of system \eqref{seq1}-\eqref{sN}.  Now, the function $h$ is defined in \eqref{h}.

\begin{prop} \label{prop:stat}
Assume that $a\leq d_c$, $\gamma >0$, and $d,d_b,d_g, \kappa_0$ be arbitrary and positive. 
Then, the only nonnegative solution of  system \eqref{seq1}-\eqref{sN}
is the trivial steady state $(0,0, \kappa_0/d_g)$.
\end{prop}

\begin{proof}
If $a<d_c$, the nonexistence of stationary solutions, except the trivial one, is guaranteed by  Proposition \ref{prop:a<dc}.

For $a= d_c$,  it follows from equation \eqref{seq1} that $U\equiv 0$ and,  by equation \eqref{seq2}, we obtain $V\equiv 0$. 
Consequently, using \eqref{seq3} we obtain the 
following boundary value problem for the function $W=W(x)$
\begin{eqnarray*}
&&\frac{1}{\gamma} W_{xx} -d_g W+\kappa_0=0, \quad x \in (0,1)\\
&& W_x(0)=W_x(1)=0.
\end{eqnarray*}
Here, the potential energy $H(w)= -d_g w^2/2 +\kappa_0 w$ has no local minimum and the corresponding equation for trajectories 
\eqref{energy} does not describe any closed curve, except $E=H(\kappa_0/d_g)$ and  the stationary point 
$(w,z)=(\kappa_0/d_g, 0)$.
\end{proof}

In the following proposition, we discuss the nonexistence of positive solutions of  the boundary value problem \eqref{sw1}-\eqref{sw2}.

\begin{prop}\label{thm:stat-1}
Assume that $ a>d_c$ and $ \gamma>0$. Let other constants in \eqref{sw1} be nonnegative and arbitrary.

\begin{itemize}
\item[i.] If $\kappa_0^2<\Theta$ ({\it cf.} \eqref{triangle}), 
then  the boundary value problem \eqref{sw1}-\eqref{sw2} has no positive solutions.

\item[ii.]
If $\kappa_0^2=\Theta$,
then the constant function $w\equiv \kappa_0/d_g$ is the only solution
of the boundary value problem \eqref{sw1}-\eqref{sw2}.
\end{itemize}
\end{prop}

\begin{proof}
Direct  calculations (see Appendix \ref{sec:sted:kinetic}  --  a discussion around equation \eqref{h-w-q})
 show that, under 
the assumption $\kappa_0^2<\Theta$,  
the function $h$ defined in \eqref{h} satisfies
$
h(w)=H'(w) <0$
for all  $w>0$.
Hence, the  associated potential energy 
$
H,
$
(see formula \eqref{H}, below) 
is
strictly decreasing and, in consequence,
  equation \eqref{energy} does not define any closed trajectory of system \eqref{sys2}.

If $\kappa_0^2=\Theta$, the function $H$ is decreasing, as well. Here, however, system  \eqref{sys2} has 
a stationary point $(\kappa_0/(2d_g), 0)$.
\end{proof}

Hence, in view of Propositions \ref{prop:stat} and  \ref{thm:stat-1}, 
we have to assume that $ a>d_c$  and $\kappa_0^2>\Theta$ to be able to construct
nonconstant stationary solutions of the boundary value problem \eqref{sw1}-\eqref{sw2}.
It follows from direct calculations (see the beginning of Appendix)  that, under these assumptions,
the numbers $\w_\pm$ from \eqref{const-pm} are the only positive zeros of
the function $h$ defined in \eqref{h}. 
Moreover,    the corresponding potential energy (namely, $H'=h$)
\begin{equation}\label{H}
H(w)=
- d_gw^2/2 - d_b
\frac{d_c^2 (d_b+d)}{(a -d_c)^2}\; \log w
 +\kappa_0w
\end{equation}
 has a local minimum in $\w_{-}$, a local maximum
in $\w_{+}$.

To describe inhomogeneous stationary solutions, let us come back to the abstract boundary value problem
\eqref{BV1}-\eqref{BV2} and recall how to
 calculate the number $T>0$ such that $w'(0)=w'(T)=0$.
Indeed, in view of equation \eqref{energy}, every  solution $w=w(x)$ to problem \eqref{BV1}-\eqref{BV2}
 satisfies the first order differential equation
\begin{equation}
w'(x) = \pm \sqrt{2(E-H(w(x)))}.\label{eq:w}
\end{equation}
Hence, choosing the upper branch in  \eqref{eq:w} and integrating with respect to $x$ we obtain
$$
T= \int_0^T \frac{w'(x)\;dx }{\sqrt{2(E-H(w(x)))}}.
$$
Next, let us introduce $w_{1,E}<w_{2,E}$ 
such that $w_{1,E}=w(0)$ and $w_{2,E}=w(T)$, hence $H(w_{1,E})=H(w_{2,E})=E$,
see Fig. \ref{fig:2}. 
Since $w=w(x)$ is  nondecreasing for $x\in [0,T]$, by the change of variables $y=w(x)$  we obtain
the following formula for $T$ as a function of~$E$:
\begin{equation}\label{TE}
T= T(E)= \int_{w_{1,E}}^{w_{2,E}} \frac{dy }{\sqrt{2(E-H(y))}}.
\end{equation}

Our classification of all nonnegative solutions of the boundary value problem \eqref{BV1}-\eqref{BV2} 
is based on 
 detailed analysis of the integral \eqref{TE}, which properties are stated in the
following lemmas.

\begin{lemma}\label{lem:min}
Assume  $H\in C^2(0,\infty)$  has a local minimum at $\w_-$ such that $H''(\w_-)>0$.
We consider every
$E>H(\w_-)$ such that there exist $w_{1,E}$ and $w_{2,E}$ with the following properties (see Fig. \ref{fig:2}):
\begin{itemize}
\item  $H(w_{1,E})=H(w_{2,E})=E$, 
 $H'(w_{1,E})\neq 0$,  $H'(w_{2,E})\neq 0$,
\item $w_{1,E}<\w_-< w_{2,E}$, 
\item  $H(w)<E$ for all $w\in(w_{1,E}, w_{2,E})$.
\end{itemize} 
Then, for all such constants $E$ the integral $T(E)$ defined in \eqref{TE}
is convergent and depends continuously on $E$. 
Moreover,  $\lim_{E\searrow H(\w_-)} T(E) = {\pi}/{\sqrt{H''(\w_{-})}}.$
\end{lemma}

\begin{proof} 
By the assumptions on the function $H$,
 the integrand of $T(E)$  from  \eqref{TE} has singularities at $w_{1,E}$ and  $w_{1,E}$, only.
 Thus, the  integral $T(E)$ is convergent by the Taylor expansion, because 
 $H'(w_{1,E})\neq 0$ and $H'(w_{2,E})\neq 0$. To show the continuous dependence of $T(E)$
 on $E$, it suffices to apply {\it e.g.} 
the Lebesgue dominated convergence theorem.

To calculate the limit  $\lim_{E\to H(\w_{-})} T(E)$, first, 
we recall that $E=H(w_{2,E})$ and 
we consider the integral
$$
T_1(E)\equiv  \int_{\w_{-}}^{w_{2,E}} \frac{dy }{\sqrt{2(H(w_{2,E})-H(y))}}
=
 \int_{0}^{w_{2,E}- \w_{-}} \frac{dy }{\sqrt{2(H(w_{2,E})-H(y+\w_{-}))}}
.
$$
 Defining  the new parameter $s=w_{2,E}-\w_{-}$ and  the shifted
 function $\wH(w) \equiv H(y+\w_{-})$, moreover,
changing  variables we obtain 
\begin{equation}\label{T1}
\lim_{E\to H(\w_{-})} T_1(E)= 
\lim_{s\to 0}
\int_{0}^{s} \frac{dy }{\sqrt{2\left(\wH(s)-\wH(y)\right)}}=
\lim_{s\to 0}
\int_{0}^{1} \frac{dy }{\sqrt{2s^{-2}\left(\wH(s)-\wH(sy)\right)}}.
\end{equation}
\shorter{
Since $\wH''(0)=H''(\w_{-})>0$, there exists $\varepsilon \in (0,1)$ such that 
 $m=\inf_{s\in [0,\varepsilon]}\wH''(s)>0$.
Moreover, for all $y\in [0,1]$ and $s\in (0,\varepsilon)$, we have $\wH'(sy)\geq 0$.
Hence,
by the Taylor expansion,
we obtain 
\begin{equation}\label{H:est}
\frac{\wH(s)-\wH(sy)}{s^2} \geq \frac{\wH'(sy)}{s}(1-y) 
\end{equation}
for all $y\in [0,1]$ and  $s\in (0,\varepsilon)$.  
Using again the Taylor formula $\wH'(sy)=\wH'(0) +\wH''(\theta) sy/2$ with some $\theta\in (0,sy)$
and the assumption $\wH'(0)=H'(\w_{-})=0$ we rewrite inequality \eqref{H:est} as follow
\begin{equation*}\label{H:est:below}
\frac{\wH(s)-\wH(sy)}{s^2} \geq \frac{m}{2}y(1-y)
\end{equation*}
for all $y\in [0,1]$ and  $s\in (0,\varepsilon)$.
Moreover, by a direct application of the l'Hospital rule, we have  
$$
\lim_{s\to  0} \frac{\wH(s)-\wH(sy)}{s^2} = \frac{\wH''(0)}{2} (1-y^2).
$$
}
Consequently, we apply  the Lebesgue dominated convergence theorem
 combined with the l'Hospital rule %
to the integral on the right-hand side of \eqref{T1} to show
$$
\lim_{E\to H(\w_{-})} T_1(E) = \frac{1}{ \sqrt{\wH''(0)}} \int_0^1 \frac{dy}{\sqrt{1-y^2}}.
$$
Now, it suffices to recall that $\wH''(0)=H''(\w_-)$.

In a completely analogous way, one calculates the limit
$$
\lim_{E\to H(\w_{-})} \int_{w_{1,E}}^{\w_{-}} \frac{dy }{\sqrt{2(E-H(y))}}=
\frac{1}{ \sqrt{H''(\w_{-})}} \int_{-1}^0 \frac{dy}{\sqrt{1-y^2}}
$$
and the proof is complete because $ \int_{-1}^1 1/\sqrt{1-y^2}\;dy =\pi$.
\end{proof}

\begin{lemma}\label{lem:infty}
Assume that $H\in C^2(0,\infty)$ has, for certain $0<\w_-<\w_+$, the following properties
\begin{itemize}
\item  $H$ is strictly decreasing on  $(0, \w_-)$ and strictly increasing on $(\w_-, \w_+)$;

\item $H'(\w_+)=0$.

\end{itemize}
For every $E\in (H(\w_-), H(\w_+))$, choose $w_{1,E}<\w_-<w_{2,E}<\w_+$ such that
$H(w_{1,E})=H(w_{2,E})=E$, see Fig. \ref{fig:2}. Then
$
\lim_{E\nearrow H(\w_+)}  T(E)=+\infty.
$
\end{lemma}

\begin{proof}
It suffices to show that
$$
\lim_{E\nearrow H(\w_+)} 
 \int_{\w_{-}}^{w_{2,E}} \frac{dy }{\sqrt{2(H(w_{2,E})-H(y))}}=+\infty.
$$
 The function $H$ is increasing on $(\w_-,\w_+)$ and $H'(\w_+)=0$,  hence, by the Taylor expansion, we obtain
$$
H(w_{2,E})-H(y) \leq H(\w_+)-H(y)\leq \frac{\widetilde m}{2}(\w_+ - y)^2 \qquad 
\text{with}\quad  \widetilde m = \sup_{y\in (\w_-,\w_+)} |H''(y)|
$$
for all $y\in (\w_-, w_{2,E})$. Consequently,
$$
\int_{\w_{-}}^{w_{2,E}} \frac{dy }{\sqrt{2(H(w_{2,E})-H(y))}} \geq
\frac{1}{\sqrt{\widetilde m}} 
\int_{\w_{-}}^{w_{2,E}} \frac{dy }{\w_+-y} \to+\infty
\quad\text{if} \quad w_{2,E}\to \w_+. 
$$
\end{proof}

\begin{lemma}\label{lem:monotone}
Assume that a function $H$ has all properties stated in  Lemma \ref{lem:infty}.
 Suppose, moreover, that  $H\in C^3(0,\infty)$, 
$h=H'$ satisfies $h''(w)\leq 0$ for every $w\in (0, \w_{+})$,  $h''$ is not constant on any interval, and $h(\w_-)=H'(\w_-)=0$.
Then $T(E)$ defined in \eqref{TE} is a strictly increasing function of $E \in (H(\w_{-}), H(\w{+}))$. 
\end{lemma}

\begin{proof}
Obviously, it suffices to show  that $dT(E)/dE>0$ using the explicit formula \eqref{TE}. These involved calculations 
were done by Loud \cite{L59} by using a clever change of variables.
Let us recall that result in our particular case.

First, we consider the integral 
\begin{equation}\label{T1:1}
T_1(E)=  \int_{\w_{-}}^{w_{2,E}} \frac{dy }{\sqrt{2(H(w_{2,E})-H(y))}},
\end{equation}
where shifting the variable $y$  as in the proof of Lemma  \ref{lem:min}, 
  we can assume that $\w_{-}=0$, without loss of generality.
Since $H$ is strictly increasing on  $(0, \w_{+})$, denoting $H(A)=E$, it suffices to show that the following function of $A$
$$
\wT(A) \equiv T_1(H(A)) = 
\int_0^A  \frac{dy }{\sqrt{2(H(A)-H(y))}}
$$
is strictly increasing. It is calculated in \cite[Thm. 1, Eq.~(2.5)]{L59}
that 
\begin{equation}\label{DwT}
\frac{d\wT(A)}{dA} 
= -\sqrt{2} \frac{h(A)}{H(A)} \int_0^A
\left[
\frac{H(y)h'(y)}{h(y)^2}-\frac12
\right]
\frac{dy }{\sqrt{2(H(A)-H(y))}}.
\end{equation}
Now, we denote by $f(y)$ the quantity in the brackets on the right hand side of \eqref{DwT} 
and we observe that, for every $y\in [0,A]$, we have 
\begin{equation}\label{h2f}
\frac{d}{dy} \big( h^2(y) f(y)\big) =H(y)h''(y) \quad \text{and}\quad
h^2(0)f(0)=0.
\end{equation}
By the assumptions, $H(y)\ge 0$ and $h''(y) \le 0$. Therefore, it follows from \eqref{h2f} that $h^2(y)f(y) \leq 0$ and, in consequence, $f(y)\leq 0$ for all $y\in [0,A]$. 
Finally,  the right-hand side of  equation \eqref{DwT} is positive because
$h(A)\geq 0$ and $H(A)\geq 0$.

The analysis is completely analogous  in the case of 
the counterpart of the integral $T_1(E)$, where we intergrate with respect to $y\in [w_{1,E}, \w_-]$.
\end{proof}

Now, we come back to boundary-value problem \eqref{sw1}-\eqref{sw2} with the function $h$ defined in \eqref{h}, and with the 
corresponding potential energy $H$ from \eqref{H}.

\begin{theorem}\label{thm:nopattern}
Assume that  $a>d_c$, $\kappa_0^2>\Theta$, and $\gamma\in (0, \gamma_0]$, where $\gamma_0$ is defined in \eqref{gamma0}.
Then, the constant steady states $\w_\pm$ (cf.  \eqref{const-pm}) are 
the only solutions of the boundary value problem 
\eqref{sw1}-\eqref{sw2}.
\end{theorem}

\begin{proof}
Recall that every non-constant solution $w=w(x)$ of problem \eqref{sw1}-\eqref{sw2} corresponds
to a non-constant trajectory  $(w(x), z(x))$ of system \eqref{sys2} such that $z(0)=z(T)=0$ with $T=\sqrt{\gamma}$.
However, by Lemmas \ref{lem:min} and \ref{lem:monotone}, such trajectories exist only if $T>\pi/\sqrt{H''(\w_-)}= \sqrt{\gamma_0}.$  
\end{proof}

Now, we prove our main result on the existence of continuous and positive stationary solutions of
system \eqref{eq1}-\eqref{N}. 
First, we notice that, if $W=W(x)$ is a solution  of problem \eqref{sw1}-\eqref{sw2}, then so is
$\widetilde{W} (x) = W(1-x)$.
More generally, for every $A,B\in\R$, the function $\widehat W(x)\equiv W(Ax-B)$  is a solution of equation \eqref{sw1}
with $\gamma$ replaced by $A^2\gamma$, satisfying the boundary conditions $\widehat W_x(B/A)=\widehat W_x((T+B)/A)=0$.

\begin{proof}[Proof of Theorem \ref{thm:stationary}.]
The potential energy $H=H(w)$ 
associated with equation \eqref{sw1} has the form \eqref{H}
and, under the assumption $\kappa_0^2>\Theta$, it has a local minimum in $\w_{-}$, a local maximum
in $\w_{+}$, and all properties required in Lemmas  \ref{lem:min}--\ref{lem:monotone}.
Obviously, the numbers $\w_\pm$ are constant solutions of  the boundary value problem \eqref{sw1}-\eqref{sw2}.

For every $E\in (H(\w_{-}), H(\w_{+}))$,
let us consider the trajectory $(w(x), z(x))$ of system \eqref{sys2} such that $z\geq 0$
and $z(0)=z(T)=0$ (namely, the upper half of the trajectory drawn in Fig. \ref{fig:2}).
 Here,   $T=T(E)$  is defined by the integral \eqref{TE} and, by Lemmas  \ref{lem:min}--\ref{lem:monotone},
 this is a continuous and  increasing function of $E$, which takes all values from 
the half-line $(\pi/ \sqrt{H''(\w_{-})}, \infty)$.
 Due to the change of variables \eqref{change}, the function $W(x) =w(x/T)$ is the unique increasing solution 
 of the boundary value problem \eqref{sw1}-\eqref{sw2} with $\gamma=T^2$. 
On the other hand, the lower half of the trajectory drawn in Fig. \ref{fig:2},
 more precisely, the function $\widetilde W(x) =W(1-x)$, is the unique decreasing solution 
 of \eqref{sw1}-\eqref{sw2}.
 
 Now, for fixed $\gamma>\gamma_0$, let us choose the biggest $n\in\N$ such that $\gamma/n^2> \gamma_0$
and suppose that $n\geq 2$.
 For each $k\in \{2, \dots, n\}$, we consider the unique trajectory $(w_k(x), z_k(x))$
of system \eqref{sys2}  such that $z_k(x)\geq 0$ and $z_k(0)=z_k(T_k)=0$ with $T_k=\sqrt{\gamma}/k$.
Hence, $W_k(x)= w_k(x/(kT_k))$ for $x\in (0, kT_k)$ is the solution with $k$ modes 
of the boundary value problem \eqref{sw1}-\eqref{sw2}. 
Another solutions with $k$ modes is the symmetric counterpart $\widetilde W_k(x)=W(1-x)$.

These are all solutions and there exists no other solutions, because $T=T(E)$ is a continuous and 
strictly increasing function of $E$, see Lemmas  \ref{lem:min}--\ref{lem:monotone}. Hence, we  
obtain  all trajectories $(w(x), z(x))$ of system \eqref{sys2}
such that $z(0)=z(T)$
for each  $T> \pi/\sqrt{H''(\w_{-})}$.
\end{proof}

\begin{proof}[Proof of Theorem \ref{thm:disc-stationary}.]
First, we change variables as in \eqref{change}, hence, we consider weak solutions of  
system  \eqref{seq1}-\eqref{sN} (with $1/\gamma=1$) on the interval $[0,T]$  (with $T=\sqrt\gamma$).
In these new varables, 
if $U(x)=V(x)=0$ for some $x\in (0,T)$, 
the function $W$ satisfies the equation
$$
W''(x) -d_g W(x) +\kappa_0=0.
$$
Trajectories of the corresponding system on the $W\! Z$-plane 
\begin{equation}\label{sys2-0}
W'=Z, \quad Z'= d_gW -\kappa_0
\end{equation}
are unbounded and have the following explicit form $Z=\pm \sqrt{2 E_2+d_gW^2 -2 \kappa_0 W}$ for every constant $E_2>0$, see Fig. \ref{fig:3}.

On the other hand, 
if $U(x)\neq 0$ and $V(x)\neq 0$ for some $x\in (0,T)$, we use relations \eqref{s-uv} for $U(x), V(x)$ together with 
equation \eqref{sw1} for $W(x)$  as well as the coresponding   system on the $W\! Z$-plane
\eqref{sys2} with $h$ defined in \eqref{h}. Recall that, under our assumptions,  
this system has closed trajectories described by the formula \eqref{energy}, 
see Fig.  \ref{fig:2}.

Now, we construct a continuous trajectory  $(W,Z)$, which corresponds to a weak solution of  \eqref{seq1}-\eqref{sN}
in the following way.
We begin  at the $W$-axis at $x=0$ and finishing at the $W$-axis as a certain $T>0$ going
along either trajectories of system \eqref{sys2-0} or trajectories of system \eqref{sys2}.
At each point  of the intersection of two trajectories of  different types, the function $W$ is $C^1$ because 
$Z=W'$ is continuous. Fig.~\ref{fig:3} shows  examples of  such trajectories.

One can easily  check that the constructed-in-this-way function $W=W(x)$ satisfies  the integral equation \eqref{W:weak} for every 
test function $\varphi\in C^1([0,1])$.
\end{proof}

%

\nofig{
%
%
%

\begin{figure}
  \setlength{\unitlength}{1mm}
{\footnotesize
\begin{picture}(80,55)
    \put(0,0){\includegraphics[width=75mm]{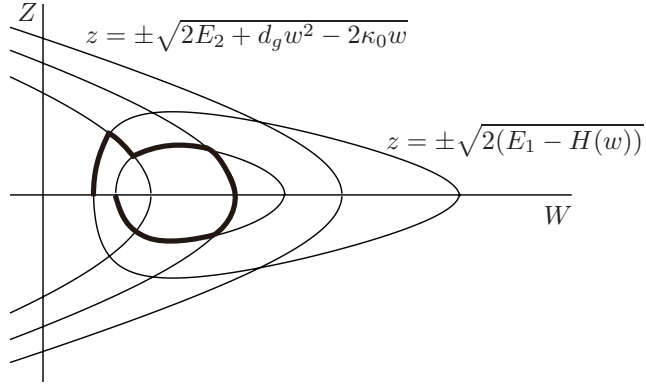}}
    \put(1,48){$Z$}
    \put(71,21){$W$}
   \put(50,31){$z=\pm \sqrt{2(E_1-H(w))}$}   
 \put(10,45){$z=\pm \sqrt{2E_2 + d_g w^2 -2\kappa_0 w}$}   
  \end{picture}
  \caption{Two closed trajectories of system \eqref{sys2} with energy levels $E_1>0$. Three trajectories of system \eqref{sys2-0}  with energy levels $E_2>0$.  The bold line shows 
 a trajectory corresponding to a discontinuous pattern.}
  \label{fig:3}
}
\end{figure}

}

\section{Linearization and spectral analysis} \label{sect:spect}

In this section, we show instability of stationary solutions of  system \eqref{eq1}-\eqref{ini}, which are constructed 
either in Theorems 
\ref{thm:stationary} or in  \ref{thm:disc-stationary}.	
  Writing equations \eqref{eq1}-\eqref{eq3} in the form
\begin{equation} \label{f123}
u_t=f_1(u,v), \quad 
v_t=f_2(u,v,w),\quad 
w_t=\frac{1}{\gamma} w_{xx}+f_3(u,v,w),
\end{equation}
we obtain that the differential of the mapping $F=(f_1,f_2,f_3):\R^3\to\R^3$ satisfies
\begin{equation}\label{DF}
DF(u,v,w)=
\left(
\begin{array}{ccc}
\frac{av^2}{(u+v)^2}-d_c   &   \frac{au^2}{(u+v)^2}  &   0  \\
2 u w                                                &         -d_b-d                     &   u^2\\
-2 u w                                               &              d                          &    -d_g-u^2
\end{array}
\right).
\end{equation}

Assume first  that $U(x)>0$ and $V(x)>0$ for some $x\in (0,1)$.
Using relations \eqref{s-uv},  we
obtain the equalities
$$
\frac{V(x)^2}{\big(U(x)+V(x)\big)^2}=\frac{d_c^2}{a^2}\qquad \text{and}\qquad 
\frac{U(x)^2}{\big(U(x)+V(x)\big)^2}=\frac{(a-d_c)^2}{a^2}.
$$
Hence, 
\begin{equation}\label{A}
DF(U,V,W)= \A (x) =(a_{ij})_{i,j=1,2,3}\equiv 
\left(
\begin{array}{ccc}
 -\frac{d_c\left( a-d_c \right)}{a}   &   \frac{(a-d_c)^2}{a}  &   0  \\
2 K                                                &         -d_b-d                     &   \frac{K^2}{W^2(x)}\\
-2 K                                              &              d                          &    -d_g- \frac{K^2}{W^2(x)}
\end{array}
\right),
\end{equation}
with the constant  
$$
K=U(x)W(x)= \frac{d_c(d_b+d)}{a-d_c},
$$
see  \eqref{s-uv}.
Notice that only the coefficients $a_{23}$ and $a_{33}$ depend on $x$.

Hence, the linearization of problem \eqref{eq1}-\eqref{ini} at a continuous steady state $(U,V,W)$ 
(namely, when $U(x)>0$ and $V(x)>0$ for all $x\in [0,1]$)
contains 
the linear operator 
\begin{equation}\label{ELL}
\L
\left(
\begin{array}{c}
 \varphi    \\
\psi   \\
\eta         
\end{array}
\right)
=
\left(
\begin{array}{ccc}
 0   &  0  &   0  \\
0     &       0       &  0\\
0         &    0     &   \frac{1}{\gamma} \partial_x^2\eta
\end{array}
\right)+\A
\left(
\begin{array}{c}
 \varphi    \\
\psi   \\
\eta         
\end{array}
\right),
\end{equation}
and we consider it 
as an operator in the Hilbert space $\H$ with the domain $D(\L)$, where
\begin{equation}\label{HDL}
\H = L^2(0,1)\times  L^2(0,1)\times L^2(0,1)\quad \text{and}\quad 
D(\L)= L^2(0,1)\times  L^2(0,1)\times W^{2,2}(0,1).
\end{equation}

\shorter{
Recall that a complex number $\lambda$ is in the {\it resolvent set}  $\rho(\L)$ of $\L$ if the range $R(\L-\lambda I)$ 
is dense in $\H$ and  $\L-\lambda I \,:\, D(\L) \to \H$ has a continuous inverse. 
The spectrum $\sigma(\L)$  (the complement of the resolvent set) 
can be  decomposed into disjoint sets: 
$\sigma (\L)= \sigma_p (\L)\cup \sigma_c (\L)\cup \sigma_r (\L)$, with 
 the {\it point spectrum}  $\sigma_p (\L)=\{\lambda\in \C\,:\, \L-\lambda I \;\; \text{is not invertible}\}$ (namely, this is the set of all eigenvalues),
 the {\it continuous spectrum}  
$\sigma_c (\L)=\{\lambda\in \C\,:\, \L-\lambda I$  has a discontinuous inverse with domain dense in  $ \H \},$
and  the  {\it residual spectrum}  
$\sigma_r (\L)=\{\lambda\in \C\,:\, \L-\lambda I$  has an inverse whose  domain is not dense in  $ \H \}$.
}


First, we prove two lemmas which allow us to characterize the point spectrum of the operator $\L$.
Here, we do not use any particular form of the coefficients in the matrix $\A$ in \eqref{A}.  Below, we denote by  $\A_{12}$
 the matrix obtained from  $\A$ after removing the third row and the third column, namely,
\begin{equation} \label{A12}
\A_{12} \equiv 
\left(
\begin{array}{cc}
a_{11}&a_{12}\\
a_{21} & a_{22}
\end{array}
\right).
\end{equation}

\begin{lemma} \label{A12:cont}
Assume that the matrix $\A_{12}$ has constant coefficients and  let
$\lambda_0$ be its  eigenvalue.
Then $\lambda_0$ belongs to the continuous spectrum of the operator $\L$ defined in \eqref{ELL}.
\end{lemma}

\begin{proof}
Let us first show that $\lambda_0$ is not an eigenvalue of the operator $\L$. 
To do it, we should show that $(\varphi, \psi, \eta)=(0,0,0)$ is the only solution of the system
\begin{equation}\label{eigen0}
\begin{array}{ccccccccc}
&&(a_{11}-\lambda_0) \varphi  & +& a_{12} \psi&    && =&0,\\
&&a_{21} \varphi & +& (a_{22} -\lambda_0) \psi &+& a_{23} \eta &=&0,\\
\frac{1}{\gamma}  \eta''  &+& a_{31} \varphi  & +&a_{32}\psi &+&(a_{33}-\lambda_0) \eta &=&0,\\
&& &&&&\eta'(0)=\eta'(1)&=&0.
\end{array}
\end{equation}
Since $\lambda_0$ is an eigenvalue of the matrix $\A_{12}$, the vectors $(a_{11}-\lambda_0, a_{12})$ and
$(a_{21},a_{22} -\lambda_0)$ are linearly dependent, namely, there is $r_1\in\R$ such that 
\begin{equation}\label{r1:vect}
(a_{11}-\lambda_0, a_{12}) =r_1 (a_{21}, a_{22} -\lambda_0).
\end{equation}
 Hence, it follows from
the first and the second equation in \eqref{eigen0} that $\eta\equiv 0$ and there exists a number $r$ 
(in fact, $r=-a_{12}/(a_{11}-\lambda_0)$)
such that $\varphi=r\psi$.
We substitute this relation into the third equation in \eqref{eigen0}  to obtain the equality
$(r a_{31} +a_{32} )\psi =0$
that implies $\psi\equiv 0$   and, consequently, $\varphi\equiv 0$.

Thus, we have proved that the operator $\L-\lambda_0 I \,:\, D(\L)\to \H$  ({\it cf.} \eqref{HDL})
is invertible.  However, its inverse   cannot be continuous, because $(\L-\lambda_0I)^{-1}$ 
is not defined on the whole space $\H$.
Indeed, using relation \eqref{r1:vect}
in the first and the second equation in  \eqref{eigen0},
we
obtain the equality $r_1 a_{23}\eta = r_1 g - f$. Hence,
\begin{equation}\label{inc:cont}
(r_1 g - f)/(r_1 a_{23})\in W^{2,2}(0,1)
\end{equation}
is a necessary condition for system \eqref{eigen0}  to have a solution $(\varphi, \psi,\eta)\in D(\L)$. 
Obviously, the condition in \eqref{inc:cont} is not satisfied for every $f,g\in L^2(0,1)$.
This completes the proof of Lemma  \ref{A12:cont}.
\end{proof}

\begin{lemma}\label{lem:eigen}
A complex number $\lambda$ is an eigenvalue of the operator $\L$ if and only if the following two conditions are 
satisfied
\begin{itemize}
\item $\lambda$ is not an eigenvalue of the matrix  $\A_{12}$,
\item the boundary value problem
\begin{equation}\label{zeta}
\begin{aligned}
&\frac{1}{\gamma} \eta''+ \frac{\det (\A-\lambda I)}{\det (\A_{12}-\lambda I)}\eta =0,& x\in (0,1)\\
& \eta'(0)=\eta'(1) =0 &
\end{aligned}
\end{equation}
has a nontrivial solution.
\end{itemize}
\end{lemma}

\begin{proof}
Assume that 
the number $\lambda\in \C$ is an eigenvalue of the operator $\L$.  Hence, the system
\begin{equation}\label{eigen}
\begin{array}{ccccccccc}
&&(a_{11}-\lambda) \varphi  & +& a_{12} \psi&    && =&0\\
&&a_{21} \varphi & +& (a_{22} -\lambda) \psi &+& a_{23} \eta &=&0\\
\frac{1}{\gamma}  \eta''  &+& a_{31} \varphi  & +&a_{32}\psi &+&(a_{33}-\lambda) \eta &=&0,
\end{array}
\end{equation}
supplemented with the boundary condition $\eta'(0)=\eta'(1)=0$, has a non-zero solution $(\varphi_0,\psi_0, \eta_0)$.
By Lemma \ref{A12:cont}, 
 $\lambda$ is not an eigenvalue 
of the matrix $\A_{12}$. Hence, $\det (\A_{12}-\lambda I)\neq 0$ and from the first and the second equation in \eqref{eigen}
we obtain 
\begin{equation} \label{phi-psi}
\varphi = \frac{a_{12}a_{23}}{\det (\A_{12}-\lambda I)} \eta 
\qquad \text{and}\qquad 
\psi = - \frac{(a_{11}-\lambda) a_{23}}{\det (\A_{12}-\lambda I)}\eta.
\end{equation}
Substituting these identities to the third equation  in \eqref{eigen} provides the relation 
\begin{equation}\label{eq:eta}
\frac{1}{\gamma} \eta'' +
\frac{a_{31}a_{12}a_{23} - a_{32} (a_{11}-\lambda) a_{23} + (a_{33}-\lambda) \det (\A_{12}-\lambda I)}
{\det(\A_{12}-\lambda I)}\eta =0,	
\end{equation}
which reduces to the equation in \eqref{zeta}, because the quantity
  $a_{31}a_{12}a_{23} - a_{32} (a_{11}-\lambda) a_{23} + (a_{33}-\lambda) \det (\A_{12}-\lambda I)$
  is the Laplace expansion of the determinant $\det(\A-\lambda I)$ with respect to its third column.

Now, let us prove the reverse implication. Assume that $\lambda\in\C$ is not an eigenvalue of $\A_{12}$ and denote by $\eta$ a non-zero
solution of problem \eqref{zeta}.  Hence, the vector  $(\varphi, \psi, \eta)$, where $\varphi$ and $\psi$ are defined in \eqref{phi-psi},  is the
eigenvector of the operator $\L$ corresponding to the eigenvalue $\lambda$.
\end{proof}

We are now in a position 
to show the existence of an infinite sequence of positive eigenvalues of the operator $\L$
and to prove Theorem  \ref{thm:unstab}. 
First, we consider the simplest case 
of the constant stationary solution $(U(x),V(x),W(x))= (\u_-, \v_-, \w_-)$. 
Here, we need a simple technical lemma.

\begin{lemma} \label{lem:pol}
Assume that a $3\times 3$-matrix $\A$ with real coefficients has all 3 eigenvalues with negative real parts.
Then, its characteristic polynomial satisfies
 $\det (\A-\lambda I)<0$  for all $\lambda\geq 0$. 
\end{lemma}
\begin{proof}
Since the polynomial $\det (\A-\lambda I)$ has real coefficients, it has three roots $\lambda_1, \lambda_2, \lambda_3\in \C$
satisfying  $\lambda_1\in \R$ and $\lambda_3= \overline{\lambda}_2$. Hence, 
$$
\det (\A-\lambda I) = - (\lambda-\lambda_1)(\lambda-\lambda_2)(\lambda-\overline{\lambda}_2)=
-(\lambda-\lambda_1)(|\lambda_2|^2-2 \lambda {\rm Re} \, \lambda_2 +\lambda^2).
$$
The  factor  $(\lambda-\lambda_1)$ 
is positive for all $\lambda\geq 0$, because $\lambda_1<0$ by the assumption.
The last factor on the right-hand side 
is positive for every $\lambda>0$,  because  ${\rm Re} \, \lambda_2<0$.
\end{proof}

\begin{theorem}[Instability of the constant steady state] \label{thm:instab:momog}
Denote by $\lambda_0$ the positive eigenvalue of the matrix $\A_{12}$. Consider the operator $\L$ from \eqref{ELL}
with the constant coefficient matrix $\A$ obtained from  the constant steady state  $(\u_-, \v_-, \w_-)$.
Then, there exists a sequence $\{\lambda_n\}_{n\in \N}$
of positive eigenvalues of the operator $\L$ that satisfy $\lambda_n\to \lambda_0$ as $n\to \infty$.
\end{theorem}

\begin{proof}
In view of Lemma \ref{lem:eigen}, we look for an infinite sequence of numbers $\lambda_n$ such the boundary value problem
\eqref{zeta} has a nontrivial solution. 
First, recall that 
the following eigenvalue value problem on the interval $[0,1]$
\begin{equation}\label{zeta:const}
\frac{1}{\gamma} \eta''+  \mu  \eta =0,\qquad 
\quad  \eta'(0)=\eta'(1) =0 
\end{equation}
has a nonzero solution $\eta(x)=\eta_n(x)=\cos(n\pi x)$ 
if  $\mu= \mu_n = n^2 \pi^2/\gamma$
for each $n\in \N$.

By Lemma  \ref{lem:pol}, the polynomial $\det (\A-\lambda I)$ is negative for all $\lambda\geq 0$, because
the steady state $(\u_-, \v_-, \w_-)$ is a stable solution of the kinetic system \eqref{heq1}-\eqref{heq3}, see 
Corollary \ref{app-cor1} in Appendix.

Since $\lambda_0$ is the only positive eigenvalue of the matrix $\A_{12}$,
the quadratic polynomial $\det (\A_{12}-\lambda I)$ has the following properties: 
$\det (\A_{12}-\lambda I)<0$ for $\lambda \in [0, \lambda_0)$ and 
$\det (\A_{12}-\lambda_0 I)=0$.
Hence, 
$
\det (\A-\lambda I)/{\det (\A_{12}-\lambda I)} > 0 
$
for  all $\lambda \in [0, \lambda_0)$ and the left-hand side limit satisfies
$$
\lim_{\lambda\to\lambda_0^-} \frac{\det (\A-\lambda I)}{\det (\A_{12}-\lambda I)} = +\infty.
$$
By continuity, we immediately obtain a sequence $\lambda_n \to \lambda_0$ such that
 $$\frac{\det (\A-\lambda_n I)}{\det (\A_{12}-\lambda_n I)} = \mu_n=\frac{n^2 \pi^2}{\gamma} $$
 for all sufficiently large $n \in N$.
\end{proof}

 More-or-less similar idea is used to find positive eigenvalues of the operator $\L$ in the case of non-homogeneous 
steady states. First, however, we recall  properties of eigenvalues $\mu\in\R$ of 
the boundary value problem
\begin{equation}\label{gen:prob}
\frac{1}{\gamma} \eta''+\mu q \eta =0, \quad 
 \eta'(0)=\eta'(1)=0,
\end{equation}
where $q\in L^\infty (0,1)$ is nonnegative.
It is well-known that problem \eqref{gen:prob} has a sequence of eigenvalues $\{\mu_n (q)\}_{n=1}^\infty$
satisfying
$$
0= \mu_0(q)< \mu_1(q) <  \mu_2(q) <  \mu_2(q) < \dots \to +\infty. 
$$
Moreover, these eigenvalues depend continuously on the potential $q=q(x)$ in the following sense.

\begin{lemma}\label{lem:eigen:cont}
Let $q_1,q_2\in L^\infty (0,1)$ be nonnegative and fixed. Denote by $\{\mu_n (q_i)\}_{n=1}^\infty$, $i=1,2$, the corresponding eigenvalues of problem \eqref{gen:prob}. Then,
there exists a constant $C>0$ such  that for all $n \in \N$, it holds
$$
|\mu_n(q_1)-\mu_n(q_2)|\leq \gamma C \|q_1-q_2\|_{L^2(0,1)}.
$$
Moreover, if $q_1(x)\leq q_2(x)$ a.e. then $\mu_n(q_2)\leq \mu_n(q_1)$ for each $n\in \N$.
\end{lemma}
\begin{proof}
Denoting by $G(x,y)$ the Green function of the operator $-d^2/dx^2$ with the Neumann boundary conditions,
we convert problem \eqref{gen:prob} into the integral equation
$$
\eta (x) =\gamma \mu(q) \int_0^1 G(x,y) q(y) \eta(y) \, dy, \quad x \in [0,1],
$$
where the right-hand side defines a compact, self-adjoint, nonnegative operator on $L^2(0,1)$.
Hence, properties of eigenvalues, stated in this lemma, result immediately from the abstract theory, see {\it eg.}
\cite[Cor. 5.6 and Thm. 5.7]{PS90}.
\end{proof}

\begin{proof}[Proof of Theorem \ref{thm:unstab}.]

By Lemma \ref{A12:cont},  the number $\lambda_0$ belongs to the continuous spectrum of the operator $\L$. 
Next, denoting
\begin{equation}\label{q:lambda}
q(x, \lambda) = \frac{\det (\A(x)-\lambda I)}{\det (\A_{12}-\lambda I)},
\end{equation}
we rewrite problem  \eqref{zeta} in the form
 \begin{equation}\label{gen:q}
\frac{1}{\gamma} \eta''(x)+ q(x, \lambda) \eta(x) =0, \quad 
 \eta'(0)=\eta'(1)=0.
\end{equation}
By Lemma \ref{lem:eigen},  it suffices to find a sequence $\lambda=\lambda_n \to \lambda_0$ such that 
problem \eqref{gen:q} has a nonzero solution. We shall proceed in a sequence of steps.

{\it Step 1.} There exists $\varepsilon>0$  such that 
for every $\lambda\in (\lambda_0-\varepsilon, \lambda_0)$, we have 
$$\min_{x\in [0,1] }q(x,\lambda)\equiv  \omega(\lambda)>0  \qquad \text{and} 
\qquad \omega(\lambda)\to\infty \quad \text{and}\quad 
\lambda\to\lambda_0.$$
To show this property of the potential $q(x,\lambda)$, we use its expanded form from equation \eqref{eq:eta}
\begin{equation}\label{q:expand}
q(x,\lambda) =a_{23}(x)
\frac{a_{31}a_{12} - a_{32} (a_{11}-\lambda)}
{\det(\A_{12}-\lambda I)}+a_{33}(x)-\lambda,
\end{equation}
where, by  \eqref{A},
\begin{equation}\label{a23a33}
a_{23}= \frac{K^2}{W^2(x)} \quad \text{and}\quad a_{33}(x)=  -d_g- \frac{K^2}{W^2(x)}, \quad
\text{with} \quad K= \frac{d_c(d_b+d)}{a-d_c},
\end{equation}
and other coefficients on the right-hand side are $x$-independent. Recall (see the proof of Theorem  \ref{thm:instab:momog})
that 
$\det(\A_{12}-\lambda I)<0$   for  $\lambda\in (0, \lambda_0)$   and 
$\det(\A_{12}-\lambda_0 I)=0$.

Next, we show that the numerator in the fraction on the right-hand side of \eqref{q:expand} satisfies
$R(\lambda_0)\equiv  a_{31}a_{12} - a_{32} (a_{11}-\lambda_0)<0$.
Indeed, using the explicit form of the coefficient of the matrix $\A_{12}$, we obtain
\begin{equation}\label{lambda_0}
\lambda_0 = \frac{1}{2}
\left( - \left(\frac{d_c(a-d_c)}{a} +d_b +d\right) +\sqrt{
\left(\frac{d_c(a-d_c)}{a} +d_b +d\right)^2+4   \frac{d_c(a-d_c)(d_b+d)}{a}
}\right)
\end{equation}
and
$$
R(\lambda_0) = -2  \frac{d_c(a-d_c)(d_b+d)}{a} + d\left(\frac{d_c(a-d_c)}{a}+\lambda_0\right).
$$
Hence, denoting $y=\big(d_c(a-d_c)\big)/a$ and $D= d_b+d$ leads to
$$
R(\lambda_0)=  \widetilde R(y)= -2Dy +\frac{d
\left( y-D+\sqrt{(y+D)^2+4Dy}\right)}{2}.
$$
Obviously, $\widetilde R(0)=0$. We leave for the reader to check that $d\widetilde R/dy<0$ for all $y>0$, which implies 
that $\widetilde R(y)<0$ for all $y>0$. This implies immediately that  $R(\lambda_0)<0$.

Hence, provided $\varepsilon >0$ is sufficiently small, for all $\lambda \in (\lambda_0-\varepsilon, \lambda_0)$, we obtain
$$\frac{R(\lambda)}{ \det(\A_{12}-\lambda I)}>0\qquad  \text{and}  
\qquad \frac{R(\lambda)}{ \det(\A_{12}-\lambda I)}\to \infty\quad  \text{as}\quad  \lambda\to\lambda_0.
$$
Since  $\min_{x\in [0,1] }a_{23}(x)>0$  (see \eqref{a23a33}) and since 
$a_{33}(x)$ is a bounded function, choosing  smaller $\varepsilon>0$  (if necessary)
we obtain 
$\min_{x\in [0,1] }q(x,\lambda)=\omega(\lambda)>0$ for all $\lambda \in (\lambda_0-\varepsilon, \lambda_0)$
and $\omega(\lambda)\to\infty$ for
$\lambda\to\lambda_0$.

{\it Step 2.} By the Sturm-Liouville theory and Step 1, for every $\lambda\in (\lambda_0-\varepsilon, \lambda_0)$,
there exists an increasing sequence 
$$
0<\mu_1(q(\cdot, \lambda))< \mu_2(q(\cdot, \lambda))< ...< \mu_n(q(\cdot, \lambda))< ...\to +\infty
$$
of eigenvalues of the problem
\begin{equation}\label{eigen:mu}
\begin{aligned}
&\frac{1}{\gamma} \eta''+ \mu  q(x,\lambda) \eta =0,& x\in (0,1)\\
& \eta'(0)=\eta'(1) =0. &
\end{aligned}
\end{equation}
Our goal is to show that there exists $\lambda_n \in (\lambda_0-\varepsilon, \lambda_0)$ 
such that $\mu_n(q(\cdot, \lambda_n))=1$.
Then, the corresponding eigenfunction of \eqref{eigen:mu} will be a non-zero solution of \eqref{gen:q}.

{\it Step 3.} For each $n\in \N$,  the quantity 
$\mu_n(q(\cdot,\lambda))$ is a continuous function of $\lambda$. 
Indeed, this is an immediate consequence of Lemma  \ref{lem:eigen:cont}, because $q(x,\lambda)$ is a continuous function 
of $\lambda\in (\lambda_0-\varepsilon, \lambda_0)$.

{\it Step 4.}
Let us show that 
\ $\mu_n(\lambda) \to 0$ when $\lambda\to\lambda_0$.
By Step 1, we have 
$\omega (\lambda) = \min_{x\in [0,1]} q(x,\lambda)\to \infty$ when $\lambda\in (\lambda_0-\varepsilon, \lambda_0)$ and $\lambda\to\lambda_0$.
Since $\omega(\lambda)$ is a positive constant for each $\lambda\in (\lambda_0-\varepsilon, \lambda_0)$, the eigenvalues $\mu$ of the problem
$$
\frac{1}{\gamma} \eta''+ \mu \omega(\lambda) \eta =0, \quad 
 \eta'(0)=\eta'(1)=0.
$$
are given explicitely $\mu_n( \omega(\lambda))= n^2\pi^2/ ( \omega(\lambda)\gamma)$.
 Since $\omega(\lambda)\leq q(x,\lambda)$, the comparison property for eigenvalues (stated in Lemma \ref{lem:eigen:cont})
implies
$$
\mu_n(q(\cdot,\lambda))\leq \mu_n( \omega(\lambda))=\frac{ n^2\pi^2}{  \omega(\lambda)\gamma}\to 0
\quad \text{as}\quad \lambda\to\lambda_0.
$$ 

{\it Conclusion.}
Since 
$
0< \mu_1(q(\cdot, \lambda))< \mu_2(q(\cdot, \lambda))< ...< \mu_n(q(\cdot, \lambda))< ...\to +\infty
$
and $\lambda_n(q(\cdot, \lambda))\to 0$ as $\lambda \to\lambda_0$, for every $n\in \N$,
it follows from the continuous dependence of $\mu_n(\lambda)$ on $\lambda$  (see Step 3)
 that there exists $\lambda_n\to \lambda_0$ 
such that $\mu_n(\lambda_n)=1$, provided $n$ is sufficiently large. 
This completes the proof of Theorem  \ref{thm:unstab}.
\end{proof}

\begin{proof}[Proof of Corollary \ref{cor:unstab}]
This is a consequence of Theorem  \ref{thm:unstab} and of a general result from \cite[Thm. 5.1.3]{Henry}
as it is explained in the paragraph following Corollary \ref{cor:unstab}.
\end{proof}

\begin{proof}[Proof of Corollary \ref{cor:disc:unstab}]
Here, the analysis is similar as in the case of continuous patterns, hence, we only emphasize the most important steps.

{\it Step 1.}
We fix  a weak stationary solution $(U_\I, V_\I, W_\I)$    
with a null set $\I\subset [0,1]$.
The Fr\'echet  derivative of the nonlinear mapping $\F: \H_\I\to \H_\I$  defined by the mappings $(f_1,f_2,f_3)$ in \eqref{f123}
  at the point $(U_\I, V_\I, W_\I)\in \H_\I$
has the form
$$
D\F(U_\I, V_\I, W_\I)
\left(
\begin{array}{c}
 \varphi    \\
\psi   \\
\eta         
\end{array}
\right)
=\A_\I(x)
\left(
\begin{array}{c}
 \varphi    \\
\psi   \\
\eta         
\end{array}
\right),
$$
where  $\A_\I(x)=\A(x)$  (see  the matrix in \eqref{A}) if $x\in [0,1]\setminus \I$ and
$$
\A_\I(x) =
\left(
\begin{array}{ccc}
0&0&0\\
0&0&0\\
0&0& -d_g
\end{array}
\right)
\qquad \text{if} \quad x\in \I.
$$
This results immediately from the definition of the Fr\'echet  derivative.

{\it Step 2.} Next, we study spectral properties  of 
the linear operator 
$$
\L_\I
\left(
\begin{array}{c}
 \varphi    \\
\psi   \\
\eta         
\end{array}
\right)
=
\left(
\begin{array}{ccc}
 0   &  0  &   0  \\
0     &       0       &  0\\
0         &    0     &   \frac{1}{\gamma} \partial_x^2\eta
\end{array}
\right)+\A_\I(x)
\left(
\begin{array}{c}
 \varphi    \\
\psi   \\
\eta         
\end{array}
\right),
$$
 in the Hilbert space $\H_\I$  (see \eqref{HI}) with the domain 
$ 
D(\L_\I)= L_\I^2(0,1)\times L_\I^2(0,1)\times W^{2,2}(0,1).
$
First, we notice that 
the counterpart of Lemma \ref{A12:cont} holds true and $\lambda_0$ (the positive eigenvalue of $\A_{12}$ ) belongs
to the continuous spectrum of $(\L_\I, D(\L_I))$.  To prove this property of $\lambda_0$, 
it suffices to follow the proof of Lemma  \ref{A12:cont}.

{\it Step 3.}
A complex number $\lambda$ is an eigenvalue of the operator $\L_\I$ if and only if the following two conditions are 
satisfied:
(i)  the number $\lambda$ is not an eigenvalue of the matrix  $\A_{12}$;
(ii)  the boundary value problem
\begin{equation}\label{zeta:I}
\begin{aligned}
&\frac{1}{\gamma} \eta''+ q_\I(\cdot, \lambda)\eta =0,& x\in (0,1)\\
& \eta'(0)=\eta'(1) =0 &
\end{aligned}
\end{equation}
has a nontrivial solution.
Here,  we denote
$$
q_\I(x, \lambda)=
\left\{
\begin{array}{ccl}
\frac{\det (\A(x)-\lambda I)}{\det (\A_{12}-\lambda I)} &\text{if} &  x\in [0,1]\setminus \I\\
-d_g-\lambda &\text{if} &  x\in  \I.
\end{array}
\right.
$$
This is the counterpart of Lemma  \ref{lem:eigen} with an almost identical proof.

{\it Step 4.} 
There exists a sequence $\{\lambda_n\}_{n\in \N}$
of positive eigenvalues of the operator $\L_\I$ that satisfy $\lambda_n\to \lambda_0$ as $n\to \infty$.
This is the statement of Theorem \ref{thm:unstab} written for discontinuous patterns and 
its proof follows almost the  same arguments. In particular, in Step 3 of the proof of Theorem  \ref{thm:unstab}, we should 
use the fact the the mapping $\lambda \mapsto q_\I(\cdot, \lambda)$ is continuous in the $L^2(0,1)$-norm ({\it cf.} Lemma \ref{lem:eigen:cont}).

{\it Step 5.} 
Since,  $-\L_\I$ is a sectorial operator, we complete the proof by   \cite[Thm. 5.1.3]{Henry}.
\end{proof}



\begin{appendix} \label{append1}

\section{Kinetic system}\label{sec:sted:kinetic}

Here, we briefly review results on the large time 
 dynamics of solutions to reaction-diffusion equations \eqref{eq1}-\eqref{eq3} supplemented with constant initial values,
namely, we consider the following system of ordinary differential equations  (so-called the {\it kinetic system}):
\begin{align}
&\frac{d\u}{dt}= \Big(\frac{a  \v}{\u+\v} -d_c\Big) \u , \label{heq1}\\
&\frac{d\v}{dt}= -d_b \v +\u^2 \w -d \v,\label{heq2}\\
&\frac{d\w}{dt}= -d_g \w -\u^2 \w +d \v +\kappa_0.\label{heq3}
\end{align}

\subsection*{Steady states.}
It is easy to see that \eqref{heq1}--\eqref{heq3} has the trivial steady state
$
(\u_0,\v_0,\w_0)= (0, 0, {\kappa_0}/{d_g}).
$
Here,  the right-hand side of equation \eqref{heq1} is satisfied  in the limit sense, namely, when $\v\searrow 0$ and $\u\searrow 0$.
On the other hand, assuming that $a > d_c$ and 
\begin{align}
k_0^2 \ge \Theta, \quad \text{where}\quad \Theta = 4d_g d_b \frac{d_c^2
 (d_b + d)}{(a-d_c)^2} \label{app-01}
\end{align}
system \eqref{heq1}--\eqref{heq3} has two positive equilibriums  $(\u_\pm,
\v_\pm, \w_\pm)$, where
\begin{align}
\u_\pm =\frac{a -d_c}{d_c}\; \v_\pm, \quad \v_\pm =\frac{d_c^2
 (d_b+d)}{(a -d_c)^2}\; \frac{1}{\w_\pm}, \quad \w_{\pm}=
 \frac{\kappa_0\pm \sqrt{\kappa_0^2 - \Theta}}{2 d_g}. \label{app-eq02}
\end{align}
Indeed, for  $\u\neq 0$, we obtain from equation \eqref{heq1} (with $d\u/dt=0$)  that
$
\u=\v (a -d_c)/d_c.
$
Substituting this expression to equation \eqref{heq2} (with $d\v/dt=0$), we obtain 
\begin{equation}\label{h-v}
\v=\frac{d_c^2 (d_b+d)}{(a -d_c)^2}\; \frac{1}{\w}.
\end{equation}
Finally, adding equations \eqref{heq2} to equation \eqref{heq3} 
(with $d\v/dt=d\w/dt=0$)
and  using 
expression \eqref{h-v} we obtain the quadratic equation
\begin{equation}\label{h-w-q}
d_g\w^2
 -\kappa_0 \w+
d_b
\frac{d_c^2 (d_b+d)}{(a -d_c)^2}
=0.
\end{equation}
It is clear that equation \eqref{h-w-q} has two roots $\w_\pm$ if $\kappa_0^2 > \Theta$ %
and one root $\bar{w}_\pm = \kappa_0/(2d_g)$ if $\kappa_0^2 = \Theta$, see \eqref{triangle}--\eqref{const-pm}. 

\subsection*{Boundedness of solutions.}
In the following, we denote by $(\u(t), \v (t), \w(t))$  a solution of system
\eqref{heq1}--\eqref{heq3} with a  nonnegative initial datum. %

\begin{prop}\label{app-pro1}
Every solution $(\u(t), \v
 (t), \w(t))$ of system \eqref{heq1}--\eqref{heq3} corresponding to a nonnegative initial datum
exists for all $t>0$. Moreover, it is nonnegative, bounded for $t>0$, 
 and satisfies
\begin{equation}\label{lim:app}
\limsup_{t \to \infty} \u (t)
 \le \frac{a \kappa_0}{\mu d_c},\qquad
\frac{\kappa_0}{\nu} \le \liminf_{t \to \infty}\big( \v(t) + \w(t)\big) \le
 \limsup_{t \to \infty}\big(  \v(t) + \w(t)\big)
 \le \frac{\kappa_0}{\mu},
\end{equation}
where $\nu = \max\{d_b,\, d_g\} > 0$ and $\mu = \min\{d_b,\, d_g\} > 0$.
\end{prop}

\begin{proof}
The global existence of nonnegative solutions is a consequence of Theorem \ref{thm:existence2}.
To show the relations in \eqref{lim:app}, it suffices to follow the proof of  Theorem \ref{thm:inv}.
\shorter{Let us be more precise. It follows from Lemma \ref{lem:positive}
 that a solution of \eqref{heq1}--\eqref{heq3} supplemented with a nonnegative initial condition
has to be  nonnegative. 
To show  that $(\u(t), \v (t), \w(t))$ exists globally in time,
it is sufficient to notice that it is bounded for all $t
> 0$.   Hence, 
adding equations \eqref{heq2} to \eqref{heq3} we  obtain the inequality
\[
(\v + \w )_t  \le -\mu (\v + \w ) (t) + \kappa_0,
\]
which implies that  $\v (t)$ and $\w (t)$ are bounded for all $t > 0$. 
Now, if $\v (t) \le M$ for a constant $M>0$ and for all $t > 0$, then, by equation \eqref{heq1}, 
the function $\u (t)$ satisfies
the inequality 
$
 \u_t \le a M - d_c \u
$
which implies that $\u (t)$ is bounded for all $t > 0$. }
\end{proof}
\subsection*{Convergence to the trivial steady state}

It is already proven (see Corollary \ref{cor:stab:trivial} 
and Proposition  \ref{prop:a<dc})
that the trivial steady state  $(0, 0, {\kappa_0}/{d_g})$ is locally asymptotically stable as 
a solution to reaction diffusion system \eqref{eq1}-\eqref{eq3}.
In the case of  the kinetic system, we can 
describe more precisely the convergence of solutions towards the trivial steady state.

First, we consider the case $a < d_c$, where,
by Proposition \ref{prop:a<dc},   all solutions of the reaction-diffusion system
converge to the trivial steady state.  Now, we improve this result for system \eqref{heq1}-\eqref{heq3}
\begin{theorem}\label{app-th1}
If $a < d_c$, then every positive  solution $(\u(t), \v (t), \w(t))$ converges, as $t\to\infty$,  towards
$(0, 0, {\kappa_0}/{d_g})$. 
Moreover, 
\begin{itemize}
\item if $d_c < d_b + d$, then $\v(t)/\u(t) \to 0$ as $t \to \infty$;
\item if $d_c > d_b + d$ and $a > d_c - (d_b + d)$, then $\v (t)/ \u(t)
 \to d_c (d_b + d)/(d_b + d + a - d_c)$ as $t \to \infty$;
\item if $d_c > d_b + d$ and $a \le d_c - (d_b + d)$, then $\v (t)/ \u(t)
 \to \infty$ as $t \to \infty$.
\end{itemize}

\end{theorem}

Next, we deal with the case $a = d_c$, where the result is analogous.

\begin{theorem}\label{app-th2}
If $a = d_c$, then every positive  solution $(\u(t), \v (t), \w(t))$ converges,  as $t\to\infty$,  towards
$(0, 0, {\kappa_0}/{d_g})$. 
 Moreover, 
\begin{itemize}
\item if $a < d_b + d$, then $\v(t)/\u(t) \to 0$ as $t \to \infty$;
\item if $a > d_b + d$, then $\v (t)/ \u(t)
 \to (a -(d_b + d))/(d_b + d)$ as $t \to \infty$.
\end{itemize}
\end{theorem}

For  $a > d_c$, we see from
Corollary \ref{cor:stab:trivial} that the trivial steady state is locally asymptotically stable. %
This convergence can be better described in the case of solutions of the kinetic system.
Notice that if
$\kappa_0^2 < \Theta$, there is no other nonnegative  constant steady states of system  \eqref{heq1}-\eqref{heq3}. 

\begin{theorem}\label{app-th3}
Assume that $a > d_c$ and $\kappa_0^2 < \Theta$. Then every positive  solution \\$(\u(t), \v (t), \w(t))$ converges, as $t\to\infty$,  towards
$(0, 0, {\kappa_0}/{d_g})$. Moreover, 
\begin{itemize}
\item if $d_c < d_b + d$, then $\v(t)/\u(t) \to 0$ as $t \to \infty$;
\item if $d_c > d_b + d$, then $\v (t)/ \u(t)
 \to (d_c -(d_b + d))/(d_b + d + a - d_c)$ as $t \to \infty$.
\end{itemize}
\end{theorem}

If $\kappa_0^2 \ge \Theta$,  there are positive steady states defined in 
\eqref{app-eq02}. Letting an initial datum  $u(0)$ to be sufficiently small, we obtain
that a solution converges to the trivial
steady state, again.

\begin{theorem}\label{app-th4}
Assume that $a > d_c$ and $\kappa_0^2 \ge \Theta$. 
Let $(\u(t), \v (t), \w(t))$ be a solution to \eqref{heq1}-\eqref{heq3} corresponding to a positive 
initial condition $(\u(0), \v (0), \w(0))$.
If $\u (0)$ satisfies
\[
 \u (0) < \u_+ \equiv  \frac{a-d_c}{d_c}\cdot \frac{\kappa_0 - \sqrt{\kappa_0^2 - \Theta}}{2d_b},
\]
then this  solution  
converges towards  the trivial steady state. Here, $\u_+ $ is the constant defined in \eqref{app-eq02}. 
\end{theorem}

Since  ideas of the proofs Theorems \ref{app-th1}--\ref{app-th4} are more-or-less the same,
we are going to sketch below the proof of Theorem \ref{app-th3}, only.
 Introducing new variables $X = \v / \u$ and $Y = \u
\w$ and writing, for simplicity, $u$ instead of $\u$, we obtain the
following system
\begin{align}
\frac{du}{dt} &= \left(\frac{a X}{1 + X} - d_c \right)u, \label{app-eq10}\\
\frac{dX}{dt}&= -(d_b + d + a-d_c)X + \frac{a X}{1 + X} + Y, \label{app-eq11}\\
\frac{dY}{dt} &= \left(\frac{a X}{1 + X} - d_c \right)Y - d_g Y - u^2 Y + d X
 u^2 + \kappa_0 u\label{app-eq12}
\end{align}
supplemented with positive initial conditions
\[
 u (0) = u_0, \quad X(0) = \frac{\v (0)}{u (0)} = X_0, \quad Y(0) = u
 (0)\w (0) = Y_0.
\]
If  $\u (t)$, $\v (t)$, and $\w (t)$ are positive  of all $t>0$, the functions $X(t)$ and $Y(t)$ are also
positive. %
The following lemma plays an important role in the proof.
\begin{lemma}\label{app-lem2}
Let $a>d_c$  and $\kappa_0^2 < \Theta$. If $X_0 < d_c /(a-d_c)$ and $Y_0 < d_c (d_b +
 d)/(a-d_c)$, then $X(t)$
 and $Y(t)$ satisfy
\begin{align}
 0<X(t) < \frac{d_c}{a-d_c} \quad \text{and}\quad 0<Y(t) < \frac{d_c (d_b + d)}{a-d_c}
 \qquad \text{for all}\quad 0 \le t < \infty. \label{app-eq100}
\end{align}
\end{lemma}

\begin{proof}[Sketch of the proof of Lemma \ref{app-lem2}.]
Since this lemma can be shown by a method similar to that from the proof of  Lemma \ref{lem:stab}, 
we just mention  an important point, only. 
Here, we derive the following inequality
\begin{align}
\left[Y(t) - \frac{d_c (d_b + d)}{a-d_c}\right]_t \le -(d_g + u_0^2)
 \left[Y(t) - \frac{d_c (d_b + d)}{a-d_c}\right] + R(u),
\label{app-eq13}
\end{align}
where 
\begin{equation}\label{Ru}
 R(u) = - \frac{d_c d_g (d_b +
 d)}{a - d_c} - \frac{d_c d_b}{a-d_c}u^2 + \kappa_0 u. 
\end{equation}
If $R (u)\le 0$, we immediately obtain $Y(t) < d_c (d_b + d)/(a-d_c)$ for all $t >
0$.  However,  it is easy to check 
that the assumption  $\kappa_0^2 < \Theta$ implies  $R(u) < 0$ for all $u > 0$. 
\end{proof}

Now, notice that, under the assumptions  $a > d_c$ and $\kappa_0^2 < \Theta$,
the inequality $X (t) < d_c /(a-d_c)$ implies
  $u_t(t) < 0$ (see eq. \eqref{app-eq10}). Consequently, if $X_0 < d_c
 /(a-d_c)$ and $Y_0 < d_c (d_b + d)/ (a-d_c)$, Lemma \ref{app-lem2} provides an estimate of a solution %
$(\u(t), \v(t), \w(t))$  for $t>0$, which allows us to show its convergence towards the trivial steady state %
(see the  proof of Theorem \ref{thm:stab:trivial}). %
On the other hand, if either $X_0 > d_c /(a-d_c)$ or $Y_0 > d_c (d_b + d)/ (a-d_c)$, %
one can show, by simple geometric arguments involving equations
\eqref{app-eq10}--\eqref{app-eq11} and %
the boundedness of solutions to system \eqref{heq1}--\eqref{heq3}, that there exists $T > 0$ such that %
$X(T) \le d_c /(a-d_c)$ and $Y (T) \le d_c (d_b + d)/ (a-d_c)$. %
It is easy to see that every solution of \eqref{app-eq10}--\eqref{app-eq12} stays in the rectangle \eqref{app-eq100}.

In order to show the convergence rate of $\u$ and $\v$ which are stated in Theorem
\ref{app-th3}, we consider the steady states of system
\eqref{app-eq10}--\eqref{app-eq12}. %
It is clear that $(u, X, Y) = (0, 0, 0)$ is an equilibrium. %
Moreover, noting that $Y(t) \to 0$ as $t \to \infty$ when $u(t) \to 0$ as $t \to \infty$, we have the following {\it nonnegative} steady state of system \eqref{app-eq10}--\eqref{app-eq12},
\begin{align}
 (u, X, Y) = 
\left(0, \dfrac{d_c-(d_b +d)}{d_b + d +
 a-d_c}, 0 \right)\qquad  \text{if}\quad  d_c > d_b + d.
\label{app-eq14}
\end{align}
In the following, we denote by $E_1$ the right-hand side of
\eqref{app-eq14}. %
We study stability of $E_1$ by analyzing eigenvalues of
 the corresponding Jacobian matrix at the equilibrium.

\begin{lemma}\label{app-lem3}
If $d_c < d_b + d$, then the equilibrium $(0, 0, 0)$ is asymptotically stable, while it is unstable if
 $d_c > d_b + d$. %
For $d_c > d_b + d$, the stationary solution $E_1$ is asymptotically stable. %
\end{lemma}
\begin{proof}
The Jacobian matrix $J$ of the nonlinear mapping defined by the right-hand side of
 system \eqref{app-eq10}--\eqref{app-eq12} is of the form
\[
 J =
\left(
\begin{array}{ccc}
\frac{a X}{1+X}-d_c & \frac{a u}{(1+X)^2} & 0 \\
0 & -(d_b + d + a-d_c) + \frac{a}{(1+X)^2} & 1 \\
-2uY + 2dXu + \kappa_0 & \frac{a Y}{(1 + X)^2} + du^2 &
 \frac{aX}{1+X}-d_c - d_g -u^2
\end{array}
\right).
\]
Hence,  at the stationary solution $(0, 0, 0)$, it  becomes 
\[
J_{(0, 0, 0)} =
\left(
\begin{array}{ccc}
-d_c & 0 & 0 \\
0 & -(d_b + d + a-d_c) & 1 \\
\kappa_0 & 0 & -d_c - d_g
\end{array}
\right),
\]
and it is clear that all eigenvalues of $J_{(0, 0, 0)}$ have negative real parts. %
Similarly, letting $J_{E_1}$ to be the corresponding Jacobian
matrix at the steady states $E_1$, we obtain
\[
J_{E_1} =
\left(
\begin{array}{ccc}
-(d_b + d) & 0 & 0 \\
0 & \frac{(d_b + d + a-d_c)(d_b + d -d_c)}{a} & 1 \\
\kappa_0 & 0 & -(d_b + d + d_g)
\end{array}
\right).
\]
It follows from inequality $d_c > d_b + d$ that all eigenvalues of $J_{E_1}$ are negative. 
\end{proof}

Now, the assertions in Theorem \ref{app-th3} follow from Lemmas \ref{app-lem2}
and \ref{app-lem3}. 

\begin{rem}
For $\kappa_0^2 \ge \Theta$, we have $R(u) \le 0$  ({\it cf.} \eqref{Ru}) provided $u \le \u_+$. %
Therefore,  under the condition $u(0) < \u_+$,  we  obtain the same assertion as in  Lemma
\ref{app-lem2}, now however,   under the assumption  $a > d_c$ and $\kappa_0^2 \ge \Theta$. %
Here, it suffices to use the fact that $du/dt < 0$ whenever $X (t) < d_c /(a-d_c)$. %
Repeating the reasoning from the proof of Theorem \ref{thm:stab:trivial}, we
 can show Theorem \ref{app-th4}. %
\end{rem}

\subsection*{Stability of positive steady states}
Let $(\u_\pm , \v_\pm , \w_\pm)$ be the positive steady states of system
\eqref{heq1}--\eqref{heq3} given by formulas \eqref{app-eq02}. %
To  study their stability,
we consider again system
\eqref{app-eq10}--\eqref{app-eq12},  which 
for $a > d_c$ and $\kappa_0 > \Theta$, 
 has two positive steady states
\[
 (u, X, Y)_\pm = \left(\u_\pm, \frac{d_c}{a-d_c}, \frac{d_c (d_b + d)}{a-d_c}\right).
\]
\begin{theorem}\label{app-lem4}
The vector $(u, X, Y)_-$ is an asymptotically stable stationary solution of system \eqref{app-eq10}--\eqref{app-eq12}, 
and $(u, X, Y)_+$ is unstable.
\end{theorem}

\begin{proof}
First we consider Jacobian matrix $J_+$ of the right-hand side of \eqref{app-eq10}--\eqref{app-eq12}
at $(u, X, Y)_+$:
\[
 J_{+} =
\left(
\begin{array}{ccc}
0 & \frac{(a-d_c)^2}{a}\u_+ & 0 \\
0 & -\left[d_b + d + \frac{d_c (a-d_c)}{a}\right] & 1 \\[0.3cm]
\sqrt{\kappa_0^2 - \Theta} & \frac{d_c (a-d_c)(d_b+d)}{a} + d \u_+^2 &
 -d_g - \u_+^2
\end{array}
\right). 
\]
Hence, every eigenvalue $\lambda$ of $J_+$ is determined by its characteristic equation $F(\lambda) = 0$, where 
\begin{equation}\label{app-eq18}
\begin{split}
F(\lambda) &= -\lambda^3 - \left[\u_+^2 + d_b + d + d_g + \frac{d_c
 (a-d_c)}{a}\right]\lambda^2 \\ 
&\quad - \left[(d_g + \u_+^2)\left(d_b + d +
 \dfrac{d_c (a - d_c)}{a}\right) - \left(\dfrac{d_c (a - d_c)}{a} (d_b +
 d) + d \u_+^2 \right)\right]\lambda \\
&\qquad + \sqrt{\kappa_0^2 - \Theta}\dfrac{(a-d_c)^2}{a}\u_+. 
\end{split}
\end{equation}
Since $F(0) > 0$ and $F(+\infty) = -\infty$, there exists at least one
positive real eigenvalue. This implies that $(u, X, Y)_+$ is unstable. %

Next, we investigate Jacobian matrix  $J_-$ at $(u, X, Y)_-$,
\[
 J_{-} =
\left(
\begin{array}{ccc}
0 & \frac{(a-d_c)^2}{a}\u_- & 0 \\
0 & -\left[d_b + d + \frac{d_c (a-d_c)}{a}\right] & 1 \\[0.3cm]
-\sqrt{\kappa_0^2 - \Theta} & \frac{d_c (a-d_c)(d_b+d)}{a} + d \u_-^2 &
 -d_g - \u_-^2
\end{array}
\right). 
\]
Here, calculations are more involved   to show stability of  $(u, X, Y)_-$.
We consider again the characteristic equation  $G(\lambda) = 0$, where
\begin{equation}\label{app-eq16}
\begin{split}
G(\lambda) &= -\lambda^3 -\left[\u_-^2 + d_b + d + d_g + \dfrac{d_c
 (a-d_c)}{a}\right]\lambda^2 \\
&\quad - \left[(d_g + \u_-^2)\left(d_b + d +
 \dfrac{d_c (a - d_c)}{a}\right) - \left(\dfrac{d_c (a - d_c)}{a} (d_b +
 d) + d \u_-^2 \right)\right]\lambda \\
&\qquad - \sqrt{\kappa_0^2 - \Theta}\dfrac{(a-d_c)^2}{a}\u_-. 
\end{split}
\end{equation}
Now, for simplicity, we set $G(\lambda) = -\lambda^3 - A\lambda^2 - B\lambda - \sqrt{\kappa_0^2
- \Theta}\dfrac{(a-d_c)^2}{a}\u_+$. %
Since $A^2 - 3B > 0$, the equation $G^\prime (\lambda) = 0$ has two negative roots
\[
 \lambda_1 = \frac{A + \sqrt{A^2 - 3B}}{-3}, \quad \lambda_2 = \frac{A - \sqrt{A^2 - 3B}}{-3}.
\]
This implies that $G(\lambda) = 0$ has a local minimum and a local maximum in
the half plane $\{\lambda \, :\, \lambda < 0 \}$. %
Moreover, we obtain that $G(\lambda_1) < 0$ by virtue of the inequality $A^2 - 3B
> 0$. %
Therefore, noting  that $G(0) < 0$, $G(-\infty)=+\infty$ and $G(+\infty) =
-\infty$, we have two possibilities:  
$$
\text{either} \quad (i)\quad  G(\lambda_2) \ge 0 \quad   \text{or} \quad  (ii)\quad 
 G(\lambda_2) < 0.
$$
Here, $G(\lambda_2)$ is given by
\[
 G(\lambda_2) = \lambda_2^2 \left(\frac{A + 2 \sqrt{A^2 -
 3B}}{3}\right)-\sqrt{\kappa_0^2 - \Theta} \frac{(a-d_c)^2}{a}\u_+ \,.
\]
If (i) occurs, then the equation $G(\lambda) = 0$ has three negative real roots,
which implies  immediately that $(u, X, Y)_-$ is asymptotically stable. Here, we have to notice that  the case (i) occurs if
$\kappa_0^2$ is sufficiently close to $\Theta$. %

Next, suppose $G(\lambda)$ satisfies (ii), namely,  $G(\lambda_2) < 0$. %
Then, the equation $G(\lambda) = 0$ has three roots: one is real and negative, we denote it by  $-
\lambda_0$, and others are complex numbers $\mu$ and $\bar{\mu}$. 
Therefore,  using the decomposition
$G(\lambda) = -(\lambda + \lambda_0)(\lambda - \mu)(\lambda -
 \bar{\mu})$, %
we obtain from \eqref{app-eq16} that 
\[
 A = \lambda_0 - 2\Re (\mu), \quad B = |\mu | - \lambda_0 2 \Re (\mu).
\]
Now, we consider the following quantity 
\begin{align}
G(-A) = AB - \sqrt{\kappa_0^2 - \Theta}\dfrac{(a-d_c)^2}{a} \u_+ \,. \label{app-eq15}
\end{align}
Differentiating the right-hand side of \eqref{app-eq15} with respect to  $\kappa_0$,
we obtain 
\begin{align}
 \frac{d}{d \kappa_0} \left[AB - \sqrt{\kappa_0^2 -
 \Theta}\dfrac{(a-d_c)^2}{a} \u_+ \right] > 0 \quad \text{for all}\
 \kappa_0^2 > \Theta. \label{app-eq17}
\end{align}
Since $G(-A) > 0$ at $\kappa_0^2 = \Theta$,  inequality \eqref{app-eq17} implies that  $G(-A) > 0$ for all
$\kappa_0^2 > \Theta$. %
Under the assumption (ii), we notice that the inequality $G(\lambda) > 0$ 
is satisfied only when $\lambda < - \lambda_0$. %
Therefore, we obtain the inequality  $-A < -\lambda_0$, which is equivalent to  
$\lambda_0 - 2\Re (\mu) >
\lambda_0$. This implies $\Re (\mu) < 0$. %
Consequently, if  the case  (ii) occurs, then real parts of all eigenvalues of $J_-$
are negative, which implies that $(u, X, Y)_-$ is asymptotically stable. %
\end{proof}

Since $X = \v/\u$ and $Y = \u \w$,  the steady states $(u, X, Y)_-$ and $(u, X, Y)_+$
of system \eqref{app-eq10}--\eqref{app-eq12}
correspond to the steady states $(\u_- , \v_- , \w_-)$ and $(\u_+ , \v_+ ,
\w_+)$ of \eqref{heq1}--\eqref{heq3}, respectively.  
Here,  one should recall the expressions from \eqref{app-eq02} to see that  $\u_\pm = d_c (d_b +
d)/((a-d_c)\w_\pm)$, hence, 
\[
 \u_- = \frac{a-d_c}{d_c}\cdot \frac{\kappa_0 + \sqrt{\kappa_0^2 -
 \Theta}}{2d_b} \quad\text{and}\quad  u_+ = \frac{a-d_c}{d_c}\cdot \frac{\kappa_0 - \sqrt{\kappa_0^2 -
 \Theta}}{2d_b}.
\]
Let us also notice that $\u_- > \u_+$. %

Thus, we have proved that  stationary solutions of the original kinetic system \eqref{heq1}--\eqref{heq3}
have the following stability properties.

\begin{cor}\label{app-cor1}
The vector
$(\u_- , \v_- , \w_-)$ is an asymptotically stable stationary solution of 
 system \eqref{heq1}--\eqref{heq3}, while  $(\u_+ , \v_+ , \w_+)$ is unstable.
\end{cor}

\end{appendix}


\section*{Acknowledgments}
The authors wish to express their gratitude to Steffen H\"arting
for his active interest in this work,  his several helpful comments and discussions
as well as for numerical simulations.
A.~Marciniak-Czochra was supported by European Research Council Starting Grant ``Biostruct'' and Emmy Noether Programme of German Research Council (DFG). 
The work of G.~Karch was partially supported 
by the MNiSzW grant No.~N~N201 418839 and 
the Foundation for Polish Science operated within the
Innovative Economy Operational Programme 2007-2013 funded by European
Regional Development Fund (Ph.D. Programme: Mathematical
Methods in Natural Sciences).
K.~Suzuki  acknowledges  
MEXT the Grant-in-Aid for Young Scientists (B) 20740087 and the Sumitomo
Foundation 100233.

\end{document}